%% file: expos_siam.tex
\documentclass[12 pt]{article}
\usepackage{fullpage}

\input{ex_shared.txt}

\title{On second-order cone positive systems%
		\thanks{This work was supported by the ELLIIT Excellence Center and by the Swedish Research Council through the LCCC Linnaeus Center. It was also supported by the European Research Council under the ERC Advanced Grant Agreements Switchlet n.670645 and ScalableControl n.834142 as well as by DGAPA-UNAM under the grant PAPIIT RA105518 and by SSF under the grant RIT15-0091 SoPhy."}}

\author{Christian Grussler \thanks{Department of Electrical Engineering and Computer Sciences, UC Berkeley, Berkeley, CA 
		({christian.grussler@berkeley.edu})}
	\and Anders Rantzer \thanks{Department of Automatic Control, Lund University, Lund, Sweden 
		({rantzer@control.lth.se}.)}}

\begin{document}

\maketitle

\begin{abstract}
	Internal positivity offers a computationally cheap certificate for external (input-output) positivity of a linear time-invariant system. However, the drawback with this certificate lies in its realization dependency. Firstly, computing such a realization requires to find a polyhedral cone with a potentially high number of extremal generators that lifts the dimension of the state-space representation, significantly. Secondly, not all externally positive systems posses an internally positive realization. Thirdly, in many typical applications such as controller design, system identification and model order reduction, internal positivity is not preserved. 

To overcome these drawbacks, we present a tractable sufficient certificate of external positivity based on second-order cones. {This certificate does not require any special state-space realization: if it succeeds with a possibly non-minimal realization, then it will do so with any minimal realization. While there exist systems where this certificate is also necessary, we also demonstrate how to construct systems, where both second-order and polyhedral cones as well as other certificates fail. Nonetheless, in contrast to other realization independent certificates, the present one appears to be favourable in terms of applicability and conservatism. Three applications are representatively discussed to underline its potential.} We show how the certificate can be used to find externally positive approximations of nearly externally positive systems and demonstrated that this may help to reduce system identification errors. The same algorithm is used then to design state-feedback controllers that provide closed-loop external positivity, a common approach to avoid over- and undershooting of the step response. Lastly, we present modifications to generalized balanced truncation such that external positivity is preserved for those systems, where our certificate applies. 
\end{abstract}

\input{sec_intro.txt}

\input{sec_prelim.txt}

\input{sec_main.txt}

\input{sec_alternate.txt}

\input{sec_MOR.txt}

\input{sec_example.txt}

\input{sec_concl.txt}

\input{sec_appendix.txt}

\bibliographystyle{plain}
\bibliography{refopt,refpos,refkpos}
\end{document}

%% file: ex_shared.txt
\usepackage{lipsum}
\usepackage{amsfonts}
\usepackage{graphicx}
\usepackage{epstopdf}
\usepackage[numbers,sort&compress]{natbib}

\usepackage{amsopn}
\usepackage{amsthm}
\usepackage{amsmath} %
\usepackage{amssymb}  %
\usepackage{mathptmx} %

\usepackage{bbm}
\usepackage{algorithm}
\usepackage{algorithmic}
\usepackage{pgfplots}
\pgfplotsset{compat=newest}
\usepackage{tikz}
\usepackage{tkz-euclide}
\usepackage{enumerate}
\definecolor{ao(english)}{rgb}{0.0, 0.5, 0.0}
\usepackage[colorlinks, citecolor = {ao(english)}, linkcolor = {ao(english)}]{hyperref} 
\usepackage[nameinlink]{cleveref}
\usepackage{siunitx}
\usetikzlibrary{arrows,shapes,trees,calc,positioning,patterns,decorations.pathmorphing,decorations.markings,backgrounds}
\usetikzlibrary{matrix}
\usepgfplotslibrary{groupplots}

\usepgfplotslibrary{colormaps}%
\pgfplotsset{%
	colormap={custom1}{rgb255=(170,170,255), rgb255=(170,170,255)}
}%

\pgfplotsset{%
	colormap={custom2}{rgb255=(130,130,130), rgb255=(130,130,130)}
}%

\newtheorem{thm}{Theorem}
\crefname{thm}{Theorem}{Theorems}

\newtheorem{prop}{Proposition}
\crefname{prop}{Proposition}{Propositions}

\newtheorem{lem}{Lemma}
\crefname{lem}{Lemma}{Lemmas}

\newtheorem{cor}{Corollary}
\crefname{cor}{Corollary}{Corollaries}

\newtheorem{rem}{Remark}
\crefname{rem}{Remark}{Remark}

\crefname{ass}{Assumption}{Assumption}

\crefname{conj}{Conjecture}{Conjectures}

\newtheorem{defn}{Definition}
\crefname{defn}{Definition}{Definitions}

\crefname{prob}{Problem}{Problems}
\crefname{algorithm}{Algorithm}{Algorithms}
\crefname{paper}{Paper}{Papers}
\crefname{figure}{Figure}{Figures}
\crefname{section}{Section}{Sections}
\Crefname{section}{Section}{Sections}

\usepackage{dsfont}
\let\mathbb=\mathds

\newcommand{\Rnv}{\mathbb{R}^{n}_{\geq 0}}

\newcommand{\Rnn}{\mathbb{R}^{n \times n}}
\newcommand{\Rnnn}{\mathbb{R}^{n \times n}_{\geq 0}}

\newcommand{\cone}{\textnormal{cone}}
\newcommand{\conv}{\textnormal{conv}}

\newcommand{\diag}{\textnormal{diag}}
\newcommand{\sign}{\textnormal{sign}}
\newcommand{\argmax}{\operatornamewithlimits{argmax}}

\newcommand{\blkdiag}{\textnormal{blkdiag}}
\newcommand{\trace}{\textnormal{trace}}
\newcommand{\inter}[1]{\textnormal{int}(#1)}

\newcommand{\cl}[1]{\textnormal{cl}(#1)}

\newcommand{\opts}{\star}

\newcommand{\inert}[1]{\iota(#1)}
\newcommand{\transp}{\mathsf{T}}

\definecolor{ao(english)}{rgb}{0.0, 0.5, 0.0}

\colorlet{FigColor1}{blue}
\colorlet{FigColor2}{red}
\colorlet{FigColor3}{ao(english)}
\colorlet{FigColor4}{orange}
\pgfplotsset{every axis plot/.append style={line width=1.5pt}}
	\definecolor{bluebell}{rgb}{0.74, 0.83, 0.9}
\definecolor{airforceblue}{rgb}{0.36, 0.54, 0.66}

\crefformat{equation}{\textup{#2(#1)#3}}
\crefrangeformat{equation}{\textup{#3(#1)#4--#5(#2)#6}}
\crefmultiformat{equation}{\textup{#2(#1)#3}}{ and \textup{#2(#1)#3}}
{, \textup{#2(#1)#3}}{, and \textup{#2(#1)#3}}
\crefrangemultiformat{equation}{\textup{#3(#1)#4--#5(#2)#6}}%
{ and \textup{#3(#1)#4--#5(#2)#6}}{, \textup{#3(#1)#4--#5(#2)#6}}{, and \textup{#3(#1)#4--#5(#2)#6}}

\Crefformat{equation}{#2Equation~\textup{(#1)}#3}
\Crefrangeformat{equation}{Equations~\textup{#3(#1)#4--#5(#2)#6}}
\Crefmultiformat{equation}{Equations~\textup{#2(#1)#3}}{ and \textup{#2(#1)#3}}
{, \textup{#2(#1)#3}}{, and \textup{#2(#1)#3}}
\Crefrangemultiformat{equation}{Equations~\textup{#3(#1)#4--#5(#2)#6}}%
{ and \textup{#3(#1)#4--#5(#2)#6}}{, \textup{#3(#1)#4--#5(#2)#6}}{, and \textup{#3(#1)#4--#5(#2)#6}}

\crefdefaultlabelformat{#2\textup{#1}#3}

%% file: sec_intro.txt
\section{Introduction}
Since the emergence of the famous Perron-Frobenius theorem \cite{perron1907theorie,frobenius1912matrizen}, positive operators, this is, mappings that leave a cone invariant, have attracted much interest  \cite{birkhoff1957extension,stern1991invariant,berman1979nonnegative,loewy1975positive,schneider1970cross,mostajeran2018ordering,berman1989nonnegative}. For dynamical systems, the importance of cone-invariance has been early recognized by Luenberger \cite{luenberger1979introduction}, but only in the recent years received considerable attention \cite{farina2011positive,rantzer2015scalable,tanaka2011bounded,angeli2003monotone,smith2008monotone,forni2016differentialpos,grussler2018strongly,kaczorek2012positive}. Whereas on the modelling side, this interest is based on the frequently appearing large compartmental network structures, e.g., in bio-medicine, economics and data networks ~\cite{brown1980compartmental,shorten2006positive,farina2011positive,luenberger1979introduction}, also for system analysis these systems offer a simplified treatment through their dominant dynamics \cite{rantzer2015scalable,forni2016differentialpos,forni2018dominance,angeli2003monotone,smith2008monotone,mostajeran2018positive,sootla2018operator,berman1989nonnegative}. Among linear time-invariant systems
\begin{equation}
\begin{aligned}
\dot{x}(t) = Ax(t)+Bu(t),\\
y(t) = Cx(t)+Du(t),
\label{eq: state-space}
\end{aligned}
\end{equation} 
with state $x \in \mathbb{R}^n$, input $u \in \mathbb{R}^m$ and output $y \in \mathbb{R}^k$, the convex cone of externally positive systems, this is, systems that map nonnegative inputs to nonnegative outputs, are the most prominent representatives of cone-invariant systems, because many physical quantities are by definition nonnegative. For example, $u$ may represent the inflow of a substance into a chemical reactor and $y$ the concentration of the resulting product. If in addition, the state $x$ obeys the nonnegativity constraint, the system is usually referred to as internally positive \cite{luenberger1979introduction,farina2011positive,berman1989nonnegative}. {Besides physical interpretations, external positivity also arises as a desired constraint, e.g., in the tracking error or closed-loop dynamics to avoid over- and undershooting \cite{deodhare1990design,darbha2003synthesis,lin1997nonovershooting,phillips1988conditions,blachini2018aggregates}}. 

{Only for few operations, however, e.g., serial, parallel and positive feedback interconnections, it is easy to verify that external positivity is preserved.} For many other operations, this can be a difficult task: examples include negative feedback, common model order reduction techniques \cite{grussler2012symmetry}, system identification \cite{grussler2017indentification} or the interconnection with non-positive systems as for compound systems \cite{grussler2018strongly,grussler2020variation}. {Thus, in order to be able to verify and enforce external positivity, a certificate that is both computationally and theoretically tractable is highly desirable. The main goal of our investigations is to provide such a certificate and to demonstrate its capabilities in controller-design, system identification and model order reduction.}  
 
It should be noted that any such certificate can only be sufficient as the problem is generally NP-hard \cite{blondel2002presence}. In fact, for single-input-single-output (SISO) systems ($m = k = 1$), external positivity is equivalent to the state remaining within a convex cone for all nonnegative inputs and $C$ lying in the corresponding dual cone \cite{ohta1984reachability}. In other words, certifying external positivity is as difficult as finding such an invariant cone. One completely characterized approach that seeks such a cone is the determination of an invariant polyhedral cone leading to an internally positive realization \cite{anderson1996nonnegative,farina1996existence,benvenuti2004tutorial}. Unfortunately, this approach comes with some drawbacks: (i) it may require an arbitrarily large number of extremal generators \cite{farina2011positive,farina1996existence,benvenuti2004tutorial} and as not all externally positive systems omit an internally positive realization, this procedure is not guaranteed to terminate; (ii) it is largely unknown how to combine it with other objectives such as the linear matrix inequalities (LMIs) that are typically found in controller design or model order reduction. This work overcomes these drawbacks by seeking an invariant second-order (ellipsoidal) cone, instead (see~\cref{fig:ellip}). As the invariance of such cones has been comprehensively studied \cite{loewy1975positive,stern1991exponential,stern1991invariant,hildebrand2011lmi2}, we can derive a simple, tractable, certificate, which is representable by semi-definite programming (SDP) and thus is solvable with standard convex optimization software \cite{peaucelle2002user}. {In particular, we will see that the certificate only requires a minimal realization and its simplicity makes it easy to combine with the LMI literature in control. Thus making it a practical tool beyond a posteriori  certification. This is an important distinction to other certificates \cite{sen2008external,drummond2019external,darbha2003synthesis,lie2008sufficient,jayasuriya1991class,meadows1972inline}, which also work with minimal realizations, but appear to be more restrictive or less practical for system sizes where LMIs can be efficiently solved. Further, the fact that some systems only possess invariant second-order cones \cite{farina2011positive} makes our certificate also necessary.}  

{A preliminary discussion of our certificate with focus on its merits to positivity preserving model order reduction has been reported earlier in~\cite{grussler2014modified}. Due to the increased interest in second-order cone invariance and external positivity \cite{forni2018dominance,grussler2018strongly,grussler2017indentification,zheng2016projected,altafini2016minimal,ebhihara2018analysis,zheng2019minimal,blachini2018aggregates} since then, we decided to shift the focus in this work towards the certificate itself and its applications. Further, this work complements our preliminary results with the following additional highlights:}
\begin{enumerate}[(I)]
	\item {\emph{Additional analysis and generalization:} We discuss benefits and restrictiveness of the certificate and compare it with other methods \cite{sen2008external,drummond2019external,darbha2003synthesis,farina1996existence,lie2008sufficient,jayasuriya1991class}. In particular, it is shown that there exist externally positive systems, whose positivity cannot be certified by a second-order, a polyhedral cone or any of the compared methods. This provides an incentive for the search of more general invariant cones in the future. Further, we extend our modified balanced truncation method to the use of Lyapunov inequalities, as opposed to only equalities, which as for classical generalized balanced truncation \cite{sandberg2004balanced,beck1996model} yields the familiar error bounds and allows to additionally add Lyapunov-based constraints into the reduction process.}

	\item {\emph{Approximation of nearly externally positive systems in system identification:} While the preservation of external and internal positivity in model order reduction has received attention by us and others \cite{sootla2012scalable,grussler2012symmetry,reis2009positivity}, only recently the approximation of nearly positive systems with positive ones has been considered \cite{sato2020construction}. Unfortunately, as such approximations have been constructed based on internal positivity, the drawbacks that come from this certification method are inherited, which may manifest in a slow convergence \cite{sato2020construction}. Here, we propose an alternating correction procedure (additive corrections in $A$) that finds stable systems that fulfil our certificate in order to generate externally positive approximates to arbitrary systems. As demonstrated in our case study, this can be used to account for the loss of external positivity in system identification \cite{grussler2017indentification} and even improve the quality of the identified system.} 
	
	\item {\emph{Non-over- and undershooting in state-feedback control:} 
	 The avoidance of over- and undershooting in controller design is a classical and important challenge when dealing with finite capacities, e.g., to prevent overspilling in a bottling plant. While non-overshooting only requires an externally positive error tracking system, a common way to simultaneously avoid undershooting is to design an externally positive closed loop system \cite{deodhare1990design,bement2004state,darbha2003synthesis,lin1997nonovershooting,phillips1988conditions}. In this work, we will follow this track and use the aforementioned alternating correction procedure to perform closed-loop externally positive state-feedback. As this heuristic inherits the advantages of our certificate, it avoids internal positive realizations (cf. \cite{phillips1988conditions,tanaka2011bounded}), approximations through finite discretization (cf. \cite{deodhare1990design}) as well as high-order controllers \cite{darbha2003synthesis} and is based on less restrictive criteria (cf. \cite{lin1997nonovershooting}). In our case study, it is shown that constraining the closed-loop poles as a means for a sufficiently fast response time may lead to overshooting controllers, despite the fact that a feasible controller with closed-loop external positivity can be found with our heuristic.} 
\end{enumerate}

The paper is organized as follows. First, we introduce some basic notations and preliminaries on convex cones. Subsequently, we discuss cone-invariant systems, including positive systems. Then we are set to present and discuss our first main result, the SDP-formulation of our certificate. The certificate is then used in an alternating correction procedure for controller design and the approximation of non-positive system. Subsequently, we give our second main theoretical result on generalized, positivity preserving,  balanced truncation. Finally, numerical examples are presented and a conclusion is drawn. Proofs are left to the appendix.

%% file: sec_prelim.txt
\section{Preliminaries \& Background}\label{sec:pre} 
\subsection{Notations}
Throughout this paper, we use the following notations for real matrices and vectors $X=(x_{ij}) \in \mathbb{R}^{m \times n}$. The entry-wise absolute value of $X$ is given by $|X|=(|x_{ij}|)$ and the set of entry-wise \emph{nonnegative} matrices by $\mathbb{R}^{m \times n}_{\geq 0}$. For nonnegative real-valued mappings $u: \mathbb{R}_{\geq 0} \to \mathbb{R}^{m}$, we employ the same notation and write $u(t) \in \mathbb{R}^{m}_{\geq 0}$. Submatrices of $X$ are denoted by $$X_{(p:q,s:t)} := (x_{ij})_{\substack{p \leq i \leq q,\ s \leq j \leq t}} \in \mathbb{R}^{p-q+1 \times s-t+1}$$ and accordingly $X_{(p:q,:)} := X_{(p:q,1:n)}$ and $X_{(:,s:t)}:= X_{(1:m,s:t)}$. 
$I_n$ stands for the \emph{identity matrix} in $\mathbb{R}^{n \times n}$ and $e_i$ for the \emph{$i$-th canonical unit-vector} in $\mathbb{R}^n$. For the \emph{spectrum} of $X \in \mathbb{R}^{n \times n}$, we write $\sigma(X)$, whose elements $\lambda_1(X), \dots, \lambda_n(X)$, the eigenvalues of $X$, are sorted by decreasing real part $\Re(\lambda_i(X))$ and subsorted by increasing imaginary part $\Im(\lambda_i(X))$. If $X=X^\transp$, we write  $X \succ (\succeq) 0$ for $X$ being
positive (semi-)definite, i.e., $\sigma(X)\subset[0,\infty[$. We also use these notations to describe the relation between two matrices, e.g., $A \succeq B$ defines $A-B \succeq 0$. The \emph{inertia} $\inert{X} = (i_p,i_z,i_n)$ of $X$ is defined by the number of eigenvalues with positive $i_p$, zero $i_z$ and negative $i_n$ real-parts in $\sigma(X)$.

For $\mathcal{S} \subset \mathbb{R}^{m}$, we denote its \emph{interior}, \emph{boundary} and \emph{closure} by $\inter{\mathcal{S}}$, $\partial \mathcal{S}$ and $\cl{\mathcal{S}}$, respectively. Further, we write $A \mathcal{S} := \{ Ax: x \in \mathcal{S} \}$ for its \emph{image} under $A \in \mathbb{R}^{n \times m}$, $\conv(\mathcal{S})$ and $\cone(\mathcal{S})$ for its \emph{convex hull} and \emph{convex conic hull}. %
Finally, the $H_\infty$ norm of a transfer function $G(s)$ is denoted by $\|G\|_{H_\infty}$.

\subsection{Polyhedral vs. second-order cones}
In the following let $\mathcal{K} \subset \mathbb{R}^n$ be a convex cone. $\mathcal{K}$ is called \emph{solid} if $\inter{\mathcal{K}} \neq \emptyset$ and pointed if $\mathcal{K} \cap -\mathcal{K} = \{0\}$. If it closed, solid and pointed, then $\mathcal{K}$ is referred to as \emph{proper}. The corresponding \emph{dual cone} and its interior are given by  \cite{berman1979nonnegative}
\begin{subequations}
\begin{align}
\mathcal{K}^\ast &:= \{y: y^\transp x \geq 0 \ \text{for all } x \in \mathcal{K} \}. \label{eq:dual_cone}\\
\inter{\mathcal{K}^\ast} &=  \{y: y^\transp x > 0 \ \text{for all } x \in \cl{\mathcal{K}} \setminus \{0 \} \}.
\end{align}
\end{subequations}
$\mathcal{K}$ is a \emph{polyhedral} cone if 
\begin{equation}
\mathcal{K} = \mathcal{P}_{N} := N \mathbb{R}^m_{\geq 0}
\end{equation}
for some $N \in \mathbb{R}^{n \times m}$ and a \emph{second-order/ellipsoidal} cone if 
\begin{equation}
 \mathcal{K} = \{x: \|Px\|_2 \leq c^\transp x \},
\end{equation}
for some $P \in \mathbb{R}^{m \times n}$, $c \in \mathbb{R}^n$ and $\|\cdot\|_2$ denoting the Euclidean norm. By letting $K := P^\transp P - cc^\transp$, it is easy to see that every second-order cone can alternatively be represented as
\begin{equation}
 \mathcal{K} = \mathcal{K}_{K,c} := \{x: x^\transp K x \leq 0, \ c^\transp x \geq 0\},
\end{equation}
which reveals its construction by a \emph{double-cone} 
\begin{equation}
\mathcal{K}_{K} := \{x: x^\transp K x \leq 0\} = \mathcal{K}_{K,c} \cup \mathcal{K}_{K,-c} = \mathcal{K}_{K,c} \cup -\mathcal{K}_{K,c}
\end{equation}
that is separated through a hyperplane with normal $c$ (see~\cref{fig:ellip}). In this work, we are mostly interested in \emph{proper} cones $\mathcal{K}_{K,c}$, meaning that $\inert{K} = (n-1,0,1)$ and $c$ is \emph{strictly separating}, i.e., 
\begin{subequations}
	{{
	\begin{equation}
	\lbrace x : c^\transp x =  0\rbrace \cap\mathcal{K}_K  = \lbrace 0 \rbrace, \label{eq:strict_sep}
	\end{equation}
or equivalently 
\begin{equation}
\mathcal{K}_{K,c} = \{x : x^\transp K x \leq 0, \ c^\transp x > 0 \} \cup \{ 0\}. \label{eq:strict_sep_equiv}
\end{equation}}}
\end{subequations}

\begin{lem}\label{lem:dual}
	Let $\mathcal{K}_{K,p}$ be a proper second-order cone. %
	The following are equivalent:
	\begin{enumerate}
		\item  $\mathcal{K}_{K,c} = \mathcal{K}_{K,p}$ \label{item:hyper_p}
		\item $\forall \ p^\ast \in \inter{\mathcal{K}_{K,p}}: \ c \in \inter{\mathcal{K}_{K,p}^\ast} = \inter{\mathcal{K}_{K^{-1},p^\ast}}$ 
		\label{item:dual_cone}
		\item $\exists \ x \in  \mathcal{K}_{K,p}: \ c^\transp x > 0$ and $c^\transp K^{-1} c < 0$. \label{item:dual_normal}
 		\item $\exists \ x \in  \mathcal{K}_{K,p}, \ \tau \in \mathds{R}:  \ c^\transp x > 0$ and $K+ \tau cc^\transp \succ 0$. \label{item:dual_normal_SDP}
	\end{enumerate}
\end{lem}
A proof of \cref{lem:dual} is given in \cref{proof:lem:dual}. Note that for representing a proper second-order cone through a polyhedral cone $\mathcal{P}_N$, one would need $N$ to consist of infinitely many columns \cite{boyd2004convex}.

\subsection{Cone-invariance}
\begin{defn}[$A$-invariance]
	Let $\mathcal{K} \subset \mathbb{R}^n$ and $A \in \mathbb{R}^{n \times n}$. $\mathcal{K}$ is called \emph{$A$-invariant} if $A \mathcal{K} \subset \mathcal{K}$. $\mathcal{K}$ is called \emph{exponentially $A$-invariant} if $\mathcal{K}$ is $e^{At}$-invariant for all $t\geq 0$. 
\end{defn}
\begin{rem}\label{rem:pole}
	A necessary condition for the existence of a proper convex $e^{At}$-invariant cone $\mathcal{K}$ is $\lambda_1(A) \in \mathbb{R}$ \cite{ohta1984reachability,berman1989nonnegative}. 
\end{rem}	
By \cite{ohta1984reachability}, a polyhedral cone $\mathcal{P}_{N}$ is exponentially $A$-invariant with $c \in \mathcal{P}_{N}^\ast$ if and only if
\begin{align}
\exists \gamma \geq 0, \ P \in \mathbb{R}^{m \times m}_{\geq 0}: (A+\gamma I_n) N = N P, \ N^\transp c \in \mathbb{R}^{m}_{\geq 0}. \label{eq:inv_poly}
\end{align}
A similar formulation can be derived for a proper second-order cone $\mathcal{K}_{K,c}$ (see \cref{proof:lem:inv_ellip} for a proof).
\begin{lem}
	Let $A \in \mathbb{R}^{n \times n}$ and $\mathcal{K}_{K,c} \subset \mathbb{R}^{n}$. Then, {{$\mathcal{K}_{K,c}$ is proper}} and exponentially $A$-invariant if and only if 
	\begin{equation}
	\exists \ \gamma,  \tau \in \mathbb{R}: \ A^\transp K + KA + 2 \gamma K \preceq 0, \ K + \tau cc^\transp  \succ 0. \label{eq:inv_ellip}
	\end{equation}\label{lem:inv_ellip}
\end{lem}
The following result, which is proven in \cref{proof:prop:necc_second}, shows that sometimes there only exist second-order $e^{At}$-invariant cones.
\begin{lem}
	\label{prop:necc_second}
	Let $A \in \mathbb{R}^{3 \times 3}$ with $\sigma(A) = \{\alpha,  \alpha \pm i \beta \}$ where $\alpha, \beta \in \mathbb{R}$ and $\beta \neq 0$. Then, $\mathcal{K}$ is proper, convex $e^{At}$-invariant cone if and only if $\mathcal{K} = \mathcal{K}_{K,c}$ for some $c \in \mathbb{R}^3$ and $K \in \mathbb{R}^{3 \times 3}$ with $\inert{K} = (2,0,1)$.
\end{lem}

\begin{rem} \label{rem:inertia}Assuming that $\lambda_1(A) \neq \Re(\lambda_2(A))$, there exists both, $e^{At}$-invariant polyhedral \cite{farina1996existence} and second-order cones. In fact, if $(A, A^\transp K + KA + 2 \gamma K)$ is controllable, e.g., by requiring strictness in \cref{eq:inv_ellip}, it follows that $\inert{K} = \inert{A+\gamma I_n}$ \cite{antoulas2005approximation}. Therefore, for given $A$ and $c$ with $\lambda_1(A) \neq \Re(\lambda_2(A))$, one only needs to solve  \cref{eq:inv_ellip} for some fixed \linebreak $\gamma \in (-\Re(\lambda_2(A)),-\lambda_1(A))$ in order to find a solution $(K,\tau)$ with desired inertia. This can be done by semi-definite programming \cite{boyd2004convex}. In contrast, solving \cref{eq:inv_poly} is significantly more involved, because even for fixed $\gamma$, the size of $N$ is a priori unknown and $N$ and $P$ are coupled in a non-convex fashion.  
\end{rem}

\subsection{Positive systems}
\label{sec:coneinv}
Next we discuss cone-invariant linear time-invariant systems \cref{eq: state-space}. For convenience, we will often refer to $(A,B,C,D)$ as a system, meaning that its transfer function $G(s) = C(sI_n -A)^{-1}B+D$ is realized by \cref{eq: state-space}. If $D =0$, we also write $(A,B,C)$. 
\begin{defn}[$(A,B)$-invariance]
	Let $\mathcal{K} \subset \mathbb{R}^n$. Then $\mathcal{K}$ is called  $(A,B)$-invariant if $B_{(:,j)} \in \mathcal{K}, 1 \leq j \leq m$ and $\mathcal{K}$ is exponentially $A$-invariant. 
\end{defn}
If $\mathcal{K}$ is a proper convex cone, then $(A,B)$-invariance is equivalent to $x(t) \in \mathcal{K}$ for $t \geq0$, if $u(t) \in \mathbb{R}^{m}_{\geq 0}$ and $x(0) \in \mathcal{K}$. %
The smallest $(A,B)$-invariant proper convex cone is given by
\begin{align}
\mathcal{K}_{r}(A,B) : = \cl{\cone  { \bigcup_{j=1}^m  \{e^{At} B_{(:,j)}: t\geq 0\} }},
\end{align}
the so-called \emph{reachable cone} \cite{ohta1984reachability}. One of the most frequently appearing classes of systems with $(A,B)$-invariant proper convex cones are externally and internally positive systems. 
\begin{defn}[External Positivity]
	A linear time-invariant system \cref{eq: state-space} is called \emph{externally positive} if for $x(0)=0$ all nonnegative inputs yield nonnegative outputs. 
\end{defn}
External positivity can be characterized as follows \cite{farina2011positive,ohta1984reachability}.
\begin{prop}
	\label{prop:ext_pos_equiv} Let $(A,B,C,D)$ be minimal. Then, the following are equivalent:
		\begin{enumerate}
		\item $(A,B,C,D)$ is externally positive.
		\item $\forall t\geq 0: Ce^{At}B \in \mathbb{R}^{k \times m}_{\geq 0}$ and $D \in \mathbb{R}^{k \times m}_{\geq 0}$.
		\item $C_{(i,:)} \in \mathcal{K}_r(A,B)^\ast$, $1\leq i \leq k$ and $D \in \mathbb{R}^{k \times m}_{\geq 0}$.
\item {There exists a proper convex $(A,B)$-invariant cone $\mathcal{K}$ with $C_{(i,:)}^\transp \in \mathcal{K}^\ast$, $1\leq i \leq k$ and  $D \in \mathbb{R}^{k \times m}_{\geq 0}$.}
	\end{enumerate}
\end{prop}
{The last condition is important as it allows us to certify external positivity by possibly more tractable cones than $\mathcal{K}_r(A,B)$. As such, we call a system whose external positivity can be certified by $\mathcal{K}$ also \emph{$\mathcal{K}$-positive}. A particular case are so-called \emph{internally positive} systems, where one can choose $\mathcal{K} = \Rnv$, implying the following characterization \cite{luenberger1979introduction}.}
\begin{prop}
	\label{thm: intpos}
The following are equivalent:
	\begin{enumerate}
		\item $(A,B,C,D)$ is internally positive.
		\item $\exists \alpha \geq 0: A + \alpha I \in \Rnnn$ and $B,C,D$ are element-wise nonnegative.
	\end{enumerate}
\end{prop}
It can be shown as in \cite{ohta1984reachability} that if $(A,B,C,D)$ is a $\mathcal{P}_N$-positive minimal realization, then there exists an internally positive realization. The converse also holds true for SISO systems, this is, $k=m=1$.

%% file: sec_main.txt
\section{Second-order cone positivity}\label{sec:main_ext_test}
\begin{figure}
	 \centering
	\begin{tikzpicture}
	\begin{axis}[
	view ={110}{15},
	domain=-2:2,
	y domain=0:2*pi,
	xmin=-3.3,
	xmax=3.3,
	ymin=-3,
	ymax=3,
	zmax = 2.7,
	zmin = -2.7,
	axis lines=middle,
	xtick=\empty,
	ytick=\empty,
	ztick=\empty,
	zlabel = {$x_3$},
	ylabel = {$x_2$},
	xlabel = {$x_1$},
	samples=50]

	\addplot3 [surf,z buffer=sort, fill opacity = .3, colormap name = custom2, y domain = -3:3, domain = -3.3:3.3, shader = interp, draw opacity = 0] 
	({x},
	{y},
	{0.3*y});
	\coordinate (A) at (axis cs: {1}, {-2}, {-.6});
	\coordinate (D) at (axis cs: {1}, {-2-.4*1}, {-.6+.4*3.3333});
	\draw[line width = .5 pt ,->,>=stealth](A)--(D) node[xshift = 2 pt, right] {$C^\transp$};;
	
	\addplot3 [surf,z buffer=sort, fill opacity=.3, colormap name = custom1, draw opacity=0, shader = interp] 
	({x*cos(deg(y))},
	{x*sin(deg(y))},
	{x});
	\coordinate (Q) at (axis cs:1,2.8,1.7) {};
	\coordinate (Qneg) at (axis cs:1,2.9,-1.7) {};
	\node[above] at (Q) {$\mathcal{K}_{K,C^\transp}$};
	\node[above] at (Qneg) {$\mathcal{K}_{K,-C^\transp}$};
	\addplot3+[no markers,samples=500, samples y=0, domain=0:5,variable=\t, color = red, line width = .5 pt]
	({exp(-\t)*(cos(10*deg(\t)) + sin(10*deg(\t)))},{exp(-\t)*(cos(10*deg(\t))-sin(10*deg(\t)))},{2*exp(-\t)}); \label{traj}
	\coordinate (b) at (axis cs:1,1,2) {};
	\node[circle, fill = black, inner sep=0pt,minimum size=3pt] at (b) {};
	\node[above] at (b) {$B$};
	\end{axis}
	
	\end{tikzpicture}
	\caption{Illustration of \cref{thm:ex_pos_test} for a SISO system $(A,B,C)$: $e^{At}$-invariant second-order double cone $\mathcal{K}_K = \mathcal{K}_{K,C^\transp} \cup \mathcal{K}_{K,-C^\transp}$ with strictly separating hyperplane $\{x: C x \geq 0 \}$, $B \in \mathcal{K}_{K,C^\transp}$ and \ref{traj} trajectory of $e^{At}B$ for $t \geq 0$. \label{fig:ellip}}
\end{figure}
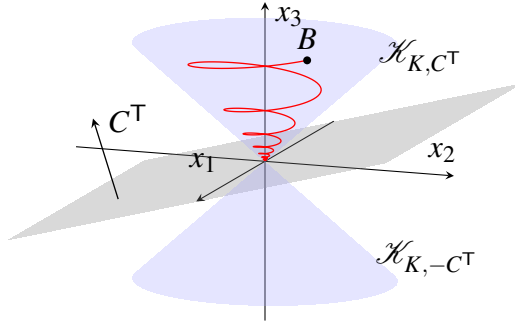

Equipped with \cref{lem:inv_ellip,prop:ext_pos_equiv}, we are ready to state our second-order cone certificate for external positivity. 
\begin{thm}[Certificate for external positivity]
	\label{thm:ex_pos_test}
	Let $(A,B,C,D)$ be a linear system and assume that there exist $K = K^\transp \in \Rnn$ and $\gamma, \tau_i \in \mathbb{R}$ such that
	\begin{subequations}
		\begin{align}
		&{A}^\transp K+ K {A} + 2\gamma K \preceq 0 \label{eq:thmtest:inv_K}\\ 
		& B_{(:,j)}^\transp K B_{(:,j)} \leq 0 \text{ for all }j\label{eq:thmtest:b_in_K}\\
		&\lambda_{n-1}(K) > 0 > \lambda_{n}(K) \label{eq:thmtest:inertia_con}\\
		& K + \tau_i C_{(i,:)}^\transp C_{(i,:)} \succ 0 \text{ for all }i \label{eq:thmtest:c_dual}\\
		&CB, \ D \in \mathbb{R}^{k\times m}_{\geq 0}\label{eq:thmtest:pos_init}
		\end{align}	
	\end{subequations}
{Then $(A,B,C,D)$ is $\mathcal{K}_{K,C_{(1,:)}}$-positive and thus externally positive with $Ce^{At}B \in \mathbb{R}^{k\times m}_{>0}$ for all $t\geq 0$.}
\end{thm}
	{A detailed proof is given \cref{proof:thm:ex_pos_test}. The certificate may be refined by applying \cref{thm:ex_pos_test} to each subsystem $(A,B_{(:,j)},C_{(i,:)},D)$, separately. However, our applications in \Cref{sec:MOR} require a common second-order cone.} \cref{fig:ellip} illustrates \cref{thm:ex_pos_test} in case of a SISO system. 
	{{\begin{rem}\label{rem:thm_main}
The assumption that $Ce^{At}B \in \mathbb{R}^{k \times m}_{> 0}$ for all $t \geq 0$ is not a strong restriction, since the sign of a floating point number can only be decided numerically up to machine precision. In particular, this condition also allows us to remove \cref{eq:thmtest:inertia_con} under mild assumptions, e.g., $\lambda_{1}(A)$ is simple. To see this, note that if there exists a $K$ fulfilling \cref{thm:ex_pos_test} with strictness in \cref{eq:thmtest:b_in_K}, then $K$ needs to have at least one negative eigenvalue, while \cref{eq:thmtest:c_dual} prevents it from having more than one. 

The existence of such a second-order cone under the assumptions of \cref{thm:ex_pos_test} can be shown as follows. For sufficiently small $\varepsilon > 0$, the system remains $e^{-A\varepsilon}\mathcal{K}_{K,C_{(1,:)}}$-positive, because 
\begin{enumerate}[i)]
	\item $e^{At}B_{(:,j)} = e^{-A\varepsilon} e^{A(t+\varepsilon)}B_{(:,j)} \in  e^{-A\varepsilon}\mathcal{K}_{K,C_{(1,:)}}$ for all $t \geq 0$.
	\item since $(e^{-A\varepsilon}\mathcal{K}_{K,C_{(1,:)}})^\ast = e^{A^\transp \varepsilon} \mathcal{K}_{K,C_{(1,:)}}^\ast$ and $C_{(i,:)}^\transp \in \inter{\mathcal{K}_{K,C_{(1,:)}}^\ast}$ (see~\cref{lem:dual} and its proof) also $C_{(i,:)}^\transp \in \inter{(e^{-A\varepsilon}\mathcal{K}_{K,C_{(1,:)}})^\ast}$.
\end{enumerate}
Further, under the assumption that $\lambda_{1}(A)$ is a simple dominant pole, $B_{(:,j)} \in \inter{e^{-A\varepsilon} \mathcal{K}_{K,C_{(1,:)}}}$. Thus, $e^{-A\varepsilon}\mathcal{K}_{K,C_{(1,:)}}$ is the desired second-order cone.
	\end{rem}}}

{{\subsection{Restrictiveness, Necessity \& Comparison}
	Next we want to study the restrictiveness and necessity of our certificate on which we will base our comparison to other certificates. We begin by considering $(A_\varepsilon,B,C)$ with
			\begin{equation}
			A_{\varepsilon} = \begin{pmatrix}
			\alpha+\varepsilon & 0 & 0\\
			0 & \alpha & \beta\\
			0 & -\beta & \alpha 
			\end{pmatrix}, \ \beta \neq 0, \ C = B^\transp, \ b_1^2 > b_2^2+b_3^2, \; \varepsilon \geq 0 \label{ex:expos_nointpos}
			\end{equation}
			Using \cref{lem:dual} and \cref{thm:ex_pos_test} with $-\alpha - \varepsilon \leq \gamma \leq -\alpha$, we can verify that the system is $\mathcal{K}_{K,C}$-positive with the \emph{Lorentz cone} $\mathcal{K}_{K,C} = \mathcal{K}_{K,C}^\ast =  \{x: x_1^2 \geq x_2^2 + x_3^2\}$. In fact, if $\varepsilon = 0$, then the system is \emph{exclusively} second-order cone-positive by \cref{prop:necc_second}.  This reveals its necessity and the fact that no certificate based on polyhedral cones such as internal positivity or \cite{altafini2016minimal} apply. {{The ability to directly deal with an arbitrary pole configuration lets our approach also appear less conservative than others \cite{jayasuriya1991class,lie2008sufficient,ball1994completely,drummond2019external}, which may require the construction of a system that under-approximates the impulse response and fulfils a certain configuration. Extreme cases of such under-approximations have been studied in \cite{meadows1972inline}.} }}
			
			Unfortunately, even when restricting ourselves to systems with $Ce^{At}B \in \mathbb{R}^{k \times m}_{> 0}$ for all $t \geq 0$, our certificate does not become a necessary condition as the following result shows. 
\begin{prop}
	\label{prop:sufficientonly}
	Let $A \in \mathbb{R}^{3 \times 3}$ be such that $\lambda_1(A) \in \mathbb{R}$, $\lambda_1(A) \neq \Re(\lambda_2(A))$ and \linebreak $\Im(\lambda_2(A)) \neq 0$. Then there exist $B, \ C^\transp, \ \Delta C^\transp \in \mathbb{R}^3$ such that
	\begin{enumerate}
		\item $\forall t \geq 0: Ce^{At}B \geq 0$, but the only $(A,B)$-invariant cone $\mathcal{K} \subset \{x: Cx \geq 0\}$ is $\mathcal{K} = \mathcal{K}_r(A,B)$, which is neither polyhedral nor second-order.
		\item $\forall t \geq 0: (C+\Delta C) e^{At}B > 0$, but no $(A,B)$-invariant cone $\mathcal{K} \subset \{x: (C+\Delta C)x \geq 0\}$ is second-order.
	\end{enumerate}
\end{prop}
\begin{rem}
	As pointed out in \cite{farina2011positive}, if $(C+\Delta C) e^{At}B > 0$ for all $t\geq0$, then $(A,B,C+\Delta C,D)$ has an internally positive realization. However, as a consequence of \cref{prop:sufficientonly}, the dimension of such a realization can still be made arbitrarily large by choosing $\Delta C$ sufficiently small. In particular, this also shows that even with the additional restriction to internally positive systems, our certificate remains only sufficient. %
\end{rem}
{{A proof to \cref{prop:sufficientonly} is stated in \cref{proof:prop:sufficientonly}, whose basic idea is illustrated in \cref{fig:counter_example}. Interestingly, similar types of systems also pose a problem for other certificates such as \cite{drummond2019external}: consider an externally positive system
\begin{equation}
G(s) = \frac{k_1}{s+p} + \frac{k_2 i}{s+\alpha + \beta i} - \frac{k_2 i}{s+\alpha - \beta i}, \; \alpha >  p >0,\; \beta, k_1 >0,
\end{equation} 
with $k_1 = 2 k_2 |\sin(\beta t^\ast) e^{-(\alpha -p )t^\ast}| < 2 k_2$, $t^\ast > 0$, i.e., the system is as in \cref{fig:counter_example}, but without the requirement that $CB = 0$. The certificate in \cite{drummond2019external} cannot be applied due to the requirement that $k_1 \geq 2 k_2$, which also remains invalid with sufficiently small perturbation as in the second item of \cref{prop:sufficientonly}.}}

	 \begin{figure}
	\centering

\begin{tikzpicture}
\begin{axis}[xmin=-1.6, xmax = 1.6,
xtick=\empty,
ytick=\empty, axis equal, axis lines=none]

\addplot[no markers,samples=500, samples y=0, domain=0:5,variable=\t, color = red, line width = .5 pt]
({exp(-\t)*sin(10*deg(\t))},{exp(-\t)*cos(10*deg(\t))}); 
\label{traj_2D}

\coordinate (b) at (axis cs:0,1) {};
\node[circle, fill = black, inner sep=0pt,minimum size=3pt] at (b) {};
\node[right, yshift =  3pt] at (b) {$B$};

\draw[line width = .5 pt, black] (axis cs: 0,1) --++(axis direction cs: 5*-0.5107 ,  5*-0.7006); \label{hyperplane}
\draw[line width = .5 pt, black] (axis cs: 0,1) --++(axis direction cs: 1.2*0.5107 ,  1.2*0.7006); 

\draw[line width = .5 pt,->,>=stealth] (axis cs:-0.5107,0.2994)  --++(axis direction cs:  .8*-0.8081,.8*0.5891) node[xshift = 2 pt, right] {$-C^\transp$}; 

\coordinate (c) at at (axis cs:-0.5107,0.2994) {};

\node[circle, fill = black, inner sep=0pt,minimum size=3pt] at (c) {};
\node[right, yshift =  3pt] at (c) {$e^{At^\ast} B$};

\draw[line width = .5 pt, dotted] (axis cs: -0.8081*0.1,1+0.5891*0.1) --++(axis direction cs: 6*-0.5107 ,  6*-0.7006); 
\draw[line width = .5 pt,dotted ] (axis cs: -0.8081*0.1,1+0.5891*0.1) --++(axis direction cs: 2*0.5107 ,  2*0.7006); \label{shift_hyper}
\draw[line width = .5 pt,->,>=stealth, dotted] (axis cs:-0.5107-0.8081*0.1 + 0.7*-0.5107  ,0.2994+0.5891*0.1-0.7006*0.7)  --++(axis direction cs:  .8*-0.8081,.8*0.5891) node[xshift = 2 pt, right] {$-(C+\Delta C)^\transp$}; 

\addplot[no markers,samples=500, samples y=0, domain=0:2*pi, variable=\t, dashed, line width = .5 pt, rotate around={54:(axis cs:0,0)}]
({2.5*sin(deg(\t))},{.67*cos(deg(\t))}); \label{thinellipse}

\end{axis}

\end{tikzpicture}
	\caption{{Illustration to \cref{prop:sufficientonly}: Projective view of $-C^\transp$ defining a hyperplane that is tangent to the convex hull of the trajectory of $e^{At}B$, i.e., $\mathcal{K}_r(A,B)$, at the linearly independent points $B$ and $e^{A^\ast}B$. If $(A,B,C)$ is $\mathcal{K}$-positive, then the hyperplane must also be tangent to $\mathcal{K}$ at both of these points. For a second-order cone, however, a hyperplane can only be tangent at two points, if they are linearly dependent. Further, there exists a arbitrarily small perturbation $\Delta C$ such that $\forall t \geq 0: (C+\Delta C) e^{At}B > 0$, which makes it necessary that any second-order that contains $B$ and $e^{A^\ast}B$ to have thinly stretched level-sets (large/small ratio between the principle axis of the defining ellipse) in order to not intersect with the hyperplane defined by $(C+\Delta C)^\transp$. As the dynamics of $e^{At}b$ are spiral, such a cone may not be $(A,B)$-invariant.} \label{fig:counter_example}}
\end{figure}
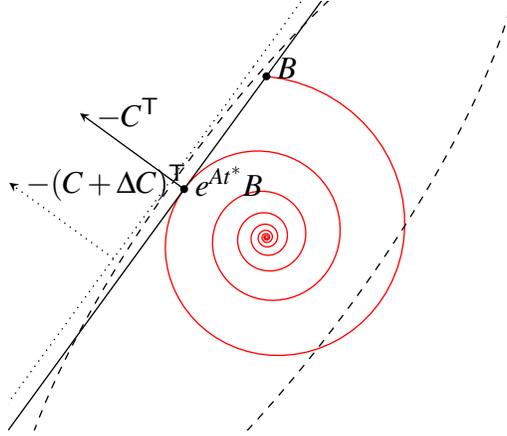

 {
	Finally, external positivity can be equivalently certified by the complete monotonicity of the transfer function \cite{sen2008external,sen2008positive,ball1994completely}. Unfortunately, while checking complete monotonicity of an impulse response is simple \cite{grussler2020variation,karlin1968total,willems1976realization}, it is a difficult task for rational functions, which makes this approach more of an analytic than implementable tool. 
	
	Overall, the indicated less restrictiveness of our certificate is bought by the need to solve LMIs. Since solving LMIs analytically may be involving, it is difficult to check whether our certificate covers any of the others completely. Further, the simplicity and analyticity of some of the other certificates is of value when it comes to large systems, where LMIs can no longer be efficiently solved. The LMI base of our certificate, however, is also an important distinction to other certificates, as it is easy to incorporate additional LMI constraints as demonstrated in the subsequent sections.  }

%% file: sec_alternate.txt
\section{Design of externally positive systems}
Designing a system such that it becomes externally positive is a desirable performance criteria. For example, a control law that results in a closed-loop externally positive systems means that monotone references signals, e.g., a step, are tracked by monotone outputs. {External positivity, thus, provides a tool to simultaneously avoid over- and undershooting, which is a highly desirable feature, e.g., in instances of limited capacities \cite{deodhare1990design,bement2004state,darbha2003synthesis,lin1997nonovershooting,phillips1988conditions}. }

Further, as external positivity is often provided through the physical quantities of our signals, it is natural to incorporate this information into modelling procedures as in system identification. Unfortunately, due to noisy measurements and other uncertainties, system identification algorithms may not produce an externally positive system \cite{grussler2017indentification}. One way of overcoming this problem is to identify a nearly externally positive system and then approximate the system with a nearby external positive one, which in turn may improve the quality of the identified model. 

In the following, we will discuss now how our certificate can help to solve these issues.  

\subsection{Alternating correction}
We start by discussing an alternating correction heuristic, which for given weight $F \in \mathbb{R}^{n \times l}$ finds a small perturbation $\Delta A \in \mathbb{R}^{l \times n}$ such that $(A+F\Delta A,B,C)$ fulfils our certificate and is asymptotically stable. In other words, we would like to solve the following non-convex problem:
\begin{equation*}
\begin{aligned}
& \underset{K,P \in \Rnn, \gamma, \tau_1, \dots, \tau_k \in \mathbb{R}}{\text{minimize}}  & & \|\Delta A\| \\
& \text{subject to} & & (A+F\Delta A)^\transp K + K(A+F\Delta A) + 2\gamma K \preceq 0\\
& &  &B_{(:,j)}^\transp K B_{(:,j)} \leq 0 \text{ for all }j \\ 
& &  & K=K^\transp, \ \lambda_{n-1}(K) > 0 > \lambda_{n}(K)\\
&  & &K + \tau_i C_{(i,:)}^\transp C_{(i,:)} \succ 0 \text{ for all }i \\
& & & (A+F\Delta A)^\transp P + P (A+F\Delta A) \preceq 0\\
& & & P \succ 0
\end{aligned}
\end{equation*}
where $\|\cdot\|$ can be any matrix norm. Next, we dualize the Lyapunov inequalities by left and right multiplication with $L = K^{-1}$ and $Q = P^{-1}$, respectively, and apply the variable changes $A_L := \Delta A L$ and $A_P := \Delta A Q$ in order to resolve these non-convex couplings. Using \cref{rem:thm_main,lem:dual} to adapt the other constraints then yields
\begin{equation*}
\begin{aligned}
& \underset{L,Q, A_L, A_Q \in \Rnn, \gamma, \tau_1, \dots, \tau_k \in \mathbb{R}}{\text{minimize}}  & & \|A_L L^{-1}\| \\
& \text{subject to} & & A L + L A^\transp + A_L^\transp F^\transp   + F A_L + 2\gamma L \preceq 0\\
& &  &C_{(j,:)}L C_{(j,:)}^\transp \leq -1 \text{ for all }j \\ 
&  & &L + \tau_i B_{(:,i)} B_{(:,i)}^\transp \succ 0 \text{ for all }i \\
& & & A Q+ Q A + A_Q^\transp F^\transp + F A_Q \preceq 0\\
& & & Q \succ 0 \\
& & & A_Q Q^{-1} = A_L L^{-1}.
\end{aligned}
\end{equation*}
Since {$\|A_L L^{-1}\| \leq \|A_L\| \|L^{-1}\|$, we may approximately keep this cost small by minimizing $\|A_L\|$. Further, the coupling of $\gamma K$ can be resolved by sweeping over different values of $\gamma$, which leaves us with the non-convexity in the last constraint.} To resolve this, we break up the problem into an alternating algorithm as outlined in \cref{alg:alt_cor}.
\begin{algorithm} 
	\caption{Find weighted externally positive stable approximation}
	\begin{algorithmic}[1]
		\STATE	\textbf{Input}: $(A,B,C)$ with $CB \in \mathbb{R}^{k\times m}_{> 0}$, $F \in \mathbb{R}^{n \times l}$, $\gamma > 0$ and precision index $\varepsilon > 0$.
		\STATE Set $e = \infty $, $\Delta A = 0$.
		\STATE \textbf{While} $e > \varepsilon$:
	\begin{equation}
		\begin{aligned}
		& \underset{L, A_L \in \Rnn, \tau_1, \dots, \tau_k \in \mathbb{R}}{\text{minimize}}  & & \|A_L\| \\
		& \text{subject to} & & A L + L A^\transp + A_L^\transp F^\transp   + F A_L + 2\gamma L \preceq 0\\
		& &  &C_{(j,:)}L C_{(j,:)}^\transp \leq -1 \text{ for all }j \\ 
		&  & &L + \tau_i B_{(:,i)} B_{(:,i)}^\transp \succ 0 \text{ for all }i \\
		\end{aligned} \label{eq:pos_opt}
		\end{equation}
		\STATE Set $\Delta A_K = A_L L^{-1}$, $\Delta A = \Delta A+ \Delta A_K$,  $A + F \Delta A_K$,
		\begin{equation}
		\begin{aligned}
		& \underset{Q, A_Q \in \Rnn}{\text{minimize}}  & & \|A_Q\| \\
		& \text{subject to}  & & A Q+ Q A + A_Q^\transp F^\transp + F A_Q \preceq 0\\
		& & & Q\succ 0
		\end{aligned} \label{eq:stab_opt}
		\end{equation}
		\STATE Set $\Delta A_P = A_Q Q^{-1}$, $\Delta A = \Delta A+ \Delta A_P$, $A = A + F \Delta A_Q$, $e = \|\Delta A_P\|$.
		\STATE \textbf{Output:} $\Delta A$
	\end{algorithmic}	\label{alg:alt_cor}
\end{algorithm}
Note that even though there is no convergence guarantee for this algorithm, in many numerical examples this procedure converges within just one iteration.

\subsection{Approximation and state-feedback controller design}
In case that $F = I_n$, \cref{alg:alt_cor} leads directly to a method of approximating a non-externally positive systems with a positive one. Since we intend to mimic the dominant dynamics of $(A,B,C)$ in this case, a reasonable range for $\gamma$ can be determined from the eigenvalues of $A$. {Although our approach has no convergence guarantees, it seems significantly less restrictive than \cite{sato2020construction}, which tries to find a minimal internally positive realization, which as a result may have slow convergence for even small dimensions.}

In case that $F = B$, our method computes a (stabilizing) state-feedback controller $u = \Delta A x(t)+r(t)$ with $r$ being a reference signal and $(A+B\Delta A,B,C)$ a closed-loop externally positive system. Here, $\gamma$ can be chosen according to the desired dominant dynamics as $-\Re(\lambda_{2}(A+B \Delta A )) > \gamma > -\lambda_1(A+B\Delta A)$. Further, \cref{eq:stab_opt} can be complemented by other LMI representable performance criteria. Additionally, our approach provides a solution to the problem of designing state-feedback controller that avoid over- and undershooting. {Other methods that accomplished solutions to this problem either rely on internal positivity \cite{phillips1988conditions,tanaka2011bounded} or other restrictive certificates \cite{lin1997nonovershooting,bement2004state}, end up with high-dimensional controllers \cite{darbha2003synthesis} or solve the problem only approximatively via finite discretization \cite{deodhare1990design}.}

%% file: sec_MOR.txt
\section{Cone balanced truncation} \label{sec:MOR}
Since many externally positive systems are formed by large networks of compartmental systems, e.g., in chemical processes and data networks \cite{luenberger1979introduction,brown1980compartmental,farina2011positive,shorten2006positive}, one often has to conduct analysis and design procedures based on their reduced order models. However, since classical model reduction techniques such as balanced truncation are not guaranteed to even preserve a dominant real pole as required in \cref{rem:pole} (unless the system is reduced to order one \cite{grussler2012symmetry}), these reduced order models do not capture some of the essential qualitative behaviours of externally positive systems such as the avoidance of over- and undershooting. Thus, potentially leading to more conservative designs as well as false conclusions. 

While \cref{alg:alt_cor} could still find a reasonable nearby externally positive approximation of priorly reduced systems, we will now introduce a modification of generalized balanced truncation (BT) \cite{beck1996model}, which even provides us with error bounds and higher quality approximations. In particular, this method can be seen as an intermediate step towards preserving internal positivity and thus retaining a compartmental structure. This is particularly important as internal positivity preserving methods \cite{reis2009positivity,sootla2012scalable} have led to rather conservative approximations \cite{grussler2012symmetry}.

We start by showing how exponential invariance with respect to a second-order cone is preserved through the concepts of cone-balanced realization and truncation. {In particular, this means that all cone-balanced truncated models preserve a dominant real pole}. Since this procedure is independent of the chosen cone, we will use it to also preserve $(A,B)$-invariance, as well as to fulfil the requirements of our external positivity certificate. Moreover, as for all generalized balanced truncation methods, it is possible to incorporate additional LMI-representable conditions, e.g., to preserve passivity. 

For a simplified exposition, we assume a minimal realization, but the readers should convince themselves that everything can be adopted to a non-minimal setting.

\begin{defn}[Cone-balanced realization]	\label{def:cone_balance}A minimal linear system realization \linebreak $(\tilde{A},\tilde{B},\tilde{C},D)$ is called \emph{cone-balanced}, if there exists diagonal $\tilde{K}$ with $\iota(\tilde{K}) = (n-1,0,1)$, diagonal $\tilde{P}, \ \tilde{Q} \succ 0$ and $\gamma \in  \mathds{R}$ such that
	\begin{subequations}
		\begin{align}
		& \tilde{A}^\transp \tilde{K} + \tilde{K} \tilde{A} + 2\gamma \tilde{K} \preceq 0, \label{eq:CBT:cone}  \\
		& \tilde{A} \tilde{P} + \tilde{P} \tilde{A}^\transp \preceq -\tilde{B}\tilde{B}^\transp, \label{eq:CBT:cont} \\
		& \tilde{A}^\transp \tilde{Q} + \tilde{Q} \tilde{A} \preceq - \tilde{C}^\transp \tilde{C},\label{eq:CBT:obs} \\
		& \tilde{p}_{11} = \tilde{q}_{11} \geq \dots \geq \tilde{p}_{nn} = \tilde{q}_{nn} \text{ and } \tilde{k}_{11} < 0. \label{eq:CBT:sort}
		\end{align}
	\end{subequations}
\end{defn}
{The idea behind \cref{def:cone_balance} is the following: By \cref{lem:dual}, $\mathcal{K}_{\tilde{K}}$ is $e^{At}$-invariant. Due to the diagonal structure of $\tilde{K}$, any element of the state, except for $x_1$, can be truncated such that the remaining systems is again exponentially $A$-invariant with respect to a second-order cone. Together with \cref{eq:CBT:cont,eq:CBT:obs}, we additionally add information on controllability and observability, which in particular allows us to provide the familiar error bound. The following result is proven in \cref{proof:thm:coneBT}.}
\begin{thm}[(Positive) cone-balanced truncation]\label{thm:coneBT}
	Suppose $(\tilde{A},\tilde{B},\tilde{C},D)$ is an asymptotically stable, cone-balanced realization of the transfer function $G(s)$ with $\tilde{K}$, $\gamma$ and $$\tilde{P} = \blkdiag \begin{pmatrix} \tilde{\sigma}_1,  \tilde{\sigma}_2 I_{l_2}\dots , \tilde{\sigma}_p I_{l_p}\end{pmatrix}$$ as in \cref{eq:CBT:cone,eq:CBT:cont},
	where $\tilde{\sigma}_2>\dots>\tilde{\sigma}_p$. 
	
	Then, for any $1 \leq r < p$, $(\tilde{A}_{(1:R:1:R)},\tilde{B}_{(1:R,:)},\tilde{C}_{(:,1:R)},D)$ with $R := 1+\sum_{i > 1}^r l_i$ and transfer function $G_R(s)$ is an asymptotically stable, cone-balanced system fulfilling \begin{align}
	\| G-G_R \|_\infty \leq 2\sum_{i =r+1}^p \tilde{\sigma}_i.  \tag{E} \label{eq:errorbound}
	\end{align}
	Further, the following are preserved:
	\begin{enumerate}
		\item $\lambda_1(\tilde{A}_{(1:R:1:R)}) \leq \gamma$. \label{item:thmCBT:gamma}
		\item  If $\mathcal{K}_{\tilde{K}}$ is $(\tilde{A},\tilde{B})$-invariant, then $\mathcal{K}_{\tilde{K}_{(1:R,1:R)}}$ is $(\tilde{A}_{(1:R:1:R)},\tilde{B}_{(1:R,:)})$-invariant \label{item:thmCBT:ABinv}
		\item If $(\tilde{A},\tilde{B},\tilde{C},D)$ is externally positive and fulfils \cref{thm:ex_pos_test} with $K = \tilde{K}$, then the same holds for $(\tilde{A}_{(1:R:1:R)},\tilde{B}_{(1:R,:)},\tilde{C}_{(:,1:R)},D)$ with $K = \tilde{K}_{(1:R,1:R)}$.
	\end{enumerate}
\end{thm}

\subsection{Cone-balancing}
Next we will discuss how to compute a cone-balanced realization. We start with the first step that yields a state-space transformation such that \cref{eq:CBT:cone} and \cref{eq:CBT:cont} are fulfilled. 
\begin{prop}\label{thm:balancing}
	Given $(A,B)$ and $N \succeq 0$, assume that there exists $\gamma > 0$, $K = K^\transp$ with $\inert{K} = (n-1,0,1)$ and $P \succ 0$ such that
	\begin{subequations}
		\begin{align}
			A^\transp K + KA + 2\gamma K \preceq 0 \label{eq:bal:cone}\\
			\trace(NK) \leq 0 \label{eq:bal:N}\\
			AP + PA^\transp =  -N \label{eq:bal:cont}
		\end{align}
	\end{subequations}
	Then there exists $T \in \mathbb{R}^{n \times n}$ such that 
	\begin{align*}
	\tilde{P} := T^{-1}PT^{-\transp} &= \blkdiag(\sigma_1,\sigma_2 I_{l_2},\dots,\sigma_s I_{l_s}) \\ 
	\tilde{K} := T^\transp KT &= \blkdiag(-\sigma_1,\sigma_2 I_{l_2},\dots,\sigma_s I_{l_s})
	\end{align*}
	where $\sigma_1 > \dots > \sigma_s > 0$, $l_2 + \dots + l_s = n - 1$ and \begin{align}
	\sigma_1^2 \geq \sum_{i > 1}l_i \sigma_i^2. \label{eq:doms1}
	\end{align}
	In particular, if $N \succeq BB^\transp$, then $(\tilde{A},\tilde{B}) = (T^{-1}AT,T^{-1}B)$ fulfills \cref{eq:CBT:cone,eq:CBT:cont} with diagonal $\tilde{K}$ and $\tilde{P}$.
\end{prop}
\cref{thm:balancing} is proven in \cref{proof:thm:balancing}. Observe that $T$ in \cref{thm:balancing} is determined in the same way as a balancing transformation in balanced truncation: the difference being that the Gramians are replaced by {{$P$ and $V|\Sigma_K|V^\transp$ with eigenvalue decomposition $K = V\Sigma_K V^\transp$}}, respectively. 
In order to decide for which states truncation causes the least error, we need another transformation of our system such that also \cref{eq:CBT:obs,eq:CBT:sort} are fulfilled. 
\begin{prop}\label{prop:delta}
	Let $(\tilde{A},\tilde{B},\tilde{C},D)$ and $\tilde{N} \succeq 0$ be such that there exist diagonal $\tilde{K}$ and $\tilde{P}$ with $\tilde{P} = |\tilde{K}|$, $\trace(\tilde{N} \tilde{K}) < 0$, $\tilde{P} \succ 0$, $\inert{K} = (n-1,0,1)$, $\mathcal{K}_{\tilde{K}}$ being $e^{\tilde{A}t}$-invariant and
	\begin{equation*}
	\tilde{A} \tilde{P} + \tilde{P} \tilde{A}^\transp \preceq - \tilde{N}.
	\end{equation*}
	Then, there exists diagonal $\Delta \succ 0$, with {diagonal entries $\delta_{11},\dots,\delta_{nn}$} such that $$\tilde{A}^\transp\Delta + \Delta \tilde{A} \preceq -\tilde{C}^\transp \tilde{C}.$$ In particular, $(\tilde{T}^{-1}\tilde{A}\tilde{T},\tilde{T}^{-1}\tilde{B},\tilde{C}\tilde{T},D)$ is cone-balanced with respect to $\mathcal{K}_{\tilde{T}^\transp\tilde{K}\tilde{T}}$, if $\tilde{N} \succeq \tilde{B}\tilde{B}^\transp$ and $\tilde{k}_{11}<0$, where
	$$\tilde{T} := \blkdiag\left(1,\frac{\tilde{p}_{22}}{\delta_{22}},\dots,\frac{\tilde{p}_{nn}}{\delta_{nn}}\right)^{\frac{1}{4}} \Pi $$
	and $\Pi$ is a permutation matrix according to \cref{eq:CBT:sort}.
\end{prop}
A proof of this result can be found in \cref{proof:prop:delta}. Since for given $K$, we can always find $N \succeq BB^\transp$ and $\tilde{N} \succeq \tilde{B}\tilde{B}^\transp$ as in \cref{thm:balancing} and \cref{prop:delta}, we have shown that \cref{eq:bal:cone} is necessary and sufficient for the existence of a cone-balanced realization. Further, if $\mathcal{K}_K$ is $(A,B)$-invariant, then 
\begin{equation}
\trace(BB^\transp K)  = \sum_j B_{(:,j)}^\transp K B_{(:,j)} \leq 0 \label{eq:BB}
\end{equation}
i.e., we can choose $N = BB^\transp$ and receive equality in \cref{eq:CBT:cont}.
\begin{cor} \label{cor:balance_contr}
	Let $(A,B,C,D)$ be {{asymptotically stable}} and $\mathcal{K}_{K}$ be $(A,B)$-invariant. Then there exists a transformation $T$ such that $(\tilde{A},\tilde{B},\tilde{C},D) := (T^{-1}AT,T^{-1}B,CT,D)$ is cone-balanced with respect to $\mathcal{K}_{T^\transp K T}$ and equality holds in \cref{eq:CBT:cont}. 
\end{cor}
Finally note that \cref{thm:balancing} and \cref{prop:delta} are intentionally presented based on $N$ and $\tilde{N}$, respectively. In this way, it is easy to see how other LMI-representable properties (see~e.g.~\cite{gugercin2004survey}) can be incorporated.

\subsection{Error-bound minimization}
\label{subsec:error}
Let us finally discuss the question of choosing $\tilde{P}$ and $\Delta$ such that the error-bound \cref{eq:errorbound} is small. We only consider the case where we also want to preserve external positivity through \cref{thm:ex_pos_test}. In this case, \cref{cor:balance_contr} applies and we can fix $\tilde{P}$ to be the controllability Gramian. Indeed, this is the best possible choice, since the eigenvalues of $\tilde{P}$ are always at least as large as those of the controllability Gramian \cite{grussler2012symmetry}. Then for finding $\Delta$, we can minimize the low-rank promoting nuclear norm \cite{fazel2001rank} of $\tilde{P}\Delta$. Alternatively, any other low-rank promoting norm \cite{grussler2016lowrank,McDonaldPS15} may also be considered. A summary of the algorithm is outlined in~\cref{alg:CBT}.

Finally note that since $\tilde{K}$ is not unique, its choice may be of considerable importance. Our experiments indicate that computing $\tilde{K}$ with respect to a balanced realization gives satisfactory results.

\begin{algorithm}
	\caption{(Positive) cone balanced truncation}
	\begin{algorithmic}[1]
		\STATE	\textbf{Input}: $(A,B,C,D)$ with $K$ that fulfills \cref{eq:thmtest:inv_K,eq:thmtest:b_in_K,eq:thmtest:inertia_con} (and \cref{eq:thmtest:c_dual,eq:thmtest:pos_init}, respectively) for (positive, respectively) cone balanced truncation.
		\STATE Find $\tilde{P}$ and $T$ in \cref{thm:balancing} with $N = BB^\transp$ and $(\tilde{A},\tilde{B},\tilde{C},D) := (T^{-1}AT,T^{-1}B,CT,D)$.
		\STATE Minimize $\sum_{i>1}{\delta_{ii}\tilde{p}_{ii}}$ subject to
		\begin{equation*}
		\tilde{A}^\transp\Delta + \Delta\tilde{A} \preceq -\tilde{C}^\transp\tilde{C}, \quad
		\Delta := \blkdiag(\delta_{11},\dots,\delta_{nn}) \succeq 0.
		\end{equation*}
		\vspace*{-.5 cm}
		\STATE Compute a cone-balanced realization $(\tilde{A},\tilde{B},\tilde{C},D)$ (see \cref{prop:delta}).
		\STATE Choose reduced order $R$ according to the predicted and desired error in \cref{eq:errorbound}.
		\STATE \textbf{Output:} $(\tilde{A}_{(1:R,1:R)},\tilde{B}_{(1:R,:)},\tilde{C}_{(:,1:R)},D)$. 
	\end{algorithmic} \label{alg:CBT}	
\end{algorithm}

%% file: sec_example.txt
\section{Case studies}\label{sec:ex}
In the following, we will illustrate our derived approaches based on case studies. We start by demonstrating that \cref{alg:alt_cor} can be used to enhance the quality of a model identified from noisy measurements. Subsequently, we design a state-feedback controller with integral action that provides an externally positive closed-loop systems. Our experiments are concluded with an example on (positive) cone balanced truncation.

\subsection{System Identification}
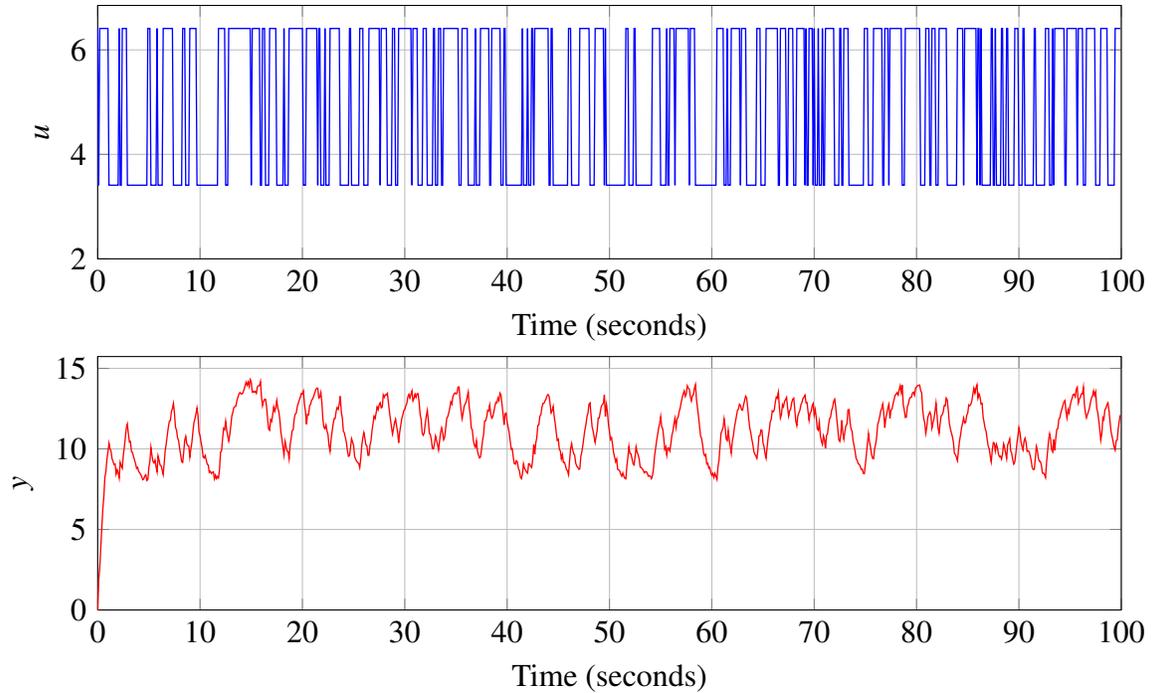
\begin{figure}
	\centering
	\begin{tikzpicture}
	\begin{groupplot}[group style={group name=my plots, group size=1 by 2, vertical sep = 1.3 cm}]
	
	\nextgroupplot[xlabel=Time (seconds),
	ylabel= $u$,
	xmin = 0,
	xmax = 100,
	grid = both,
	height = 0.3 \linewidth,
	width = .92 \linewidth,
	ymin = 2,
	scaled ticks=false, 
	]
	
	\addplot[color = FigColor1,line width = .5 pt] file{id_input.txt};
	\label{line:id_input};
		\nextgroupplot[xlabel=Time (seconds),
	ylabel= $y$,
	xmin = 0,
	xmax = 100,
	grid = both,
	height = 0.3 \linewidth,
	width = .92 \linewidth,
	ymin = 0,
	scaled ticks=false, 
	]
	
	\addplot[color = FigColor2, line width = .5 pt] file{id_output.txt};
	\label{line:id_output};

	\end{groupplot}
	\end{tikzpicture}
	\caption{Input-output data for system identifiation of the discretized heat equation model: $u$ is the same input as in the MATLAB example "dryer2"; $y$ is the corresponding output, but corrupted by white noise with variance 0.15.  \label{fig:io_id}}
\end{figure}

\begin{figure}
	\centering
	\begin{tikzpicture}
	\begin{groupplot}[group style={group name=my plots, group size=1 by 1}]
	
	\nextgroupplot[xlabel=Time (seconds),
	ylabel= Amplitude,
	xmin = 0,
	xmax = 3,
	grid = both,
	height = 0.35 \linewidth,
	width = .92 \linewidth,
	ymin = -0.13,
	ymax = 0.16,
	scaled ticks=false, 
	]
	
	\addplot[color = FigColor1,line width = 1 pt] file{impulse_id.txt};
	\label{line:id_mod};
	
	\addplot[color = FigColor2, dashdotted, line width = 1 pt] file{impulse_approx.txt};
	\label{line:ex_pos};
	
	\end{groupplot}
	\end{tikzpicture}
	\caption{Impulse response of the errors to the true system: \newline \ref{line:id_mod} N4SID of order $3$ without noise model, \newline \ref{line:ex_pos} nearby externally positive approximation to \ref{line:id_mod} through \cref{alg:alt_cor} with $\gamma = 1.25$ shows a damped behaviour, which results in an improved $H_\infty$-error of about 28 \%.   \label{fig:imp_ex_pos}}
\end{figure}
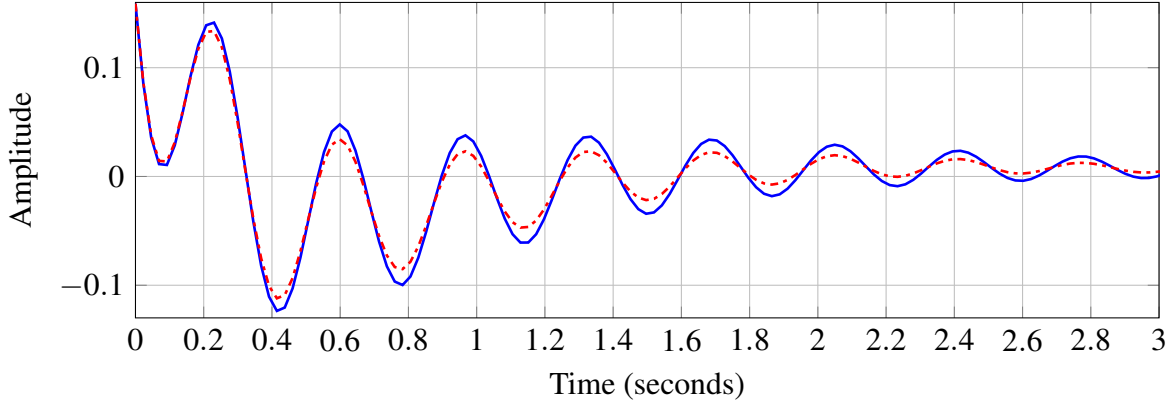

The heat equation on a two-dimensional square
\begin{equation}
\dot{T} = \triangle T = \dfrac{\partial^2}{\partial x^2}T + \dfrac{\partial^2}{\partial y^2}T \label{Heat}
\end{equation}
with control of the Dirichlet boundary conditions of the four edges yields an internally positive system if discretized on a uniform grid:
\begin{equation}
\dot{T} = A T + B u \; \text{ with $u \in \mathbb{R}^{4}$ and $T \in \mathbb{R}^{N^2}$}\; \label{HeatSS}
\end{equation}
where $A := (a_{ij}) \in \mathbb{R}^{N^2 \times N^2}$ and $B := (b_{ij}) \in \mathbb{R}^{N^2 \times 4}$ are zero except for 
\begin{align*}
a_{ii} &:= -4,  &\text{for }& \ i=1,2,\dots , N^2\\
a_{i,i+1} & = a_{i+1,i} := 1,  & \text{for }& \ i=1,\dots , N^2-1\\
a_{i,N+i} &= a_{N+i,i} := 1,  & \text{for }& \ i=1,\dots , N(N-1)
\end{align*}
and
\begin{align*}
b_{i1} &:= 1,  &\text{for }& \ i=1,2,\dots , N \\
b_{i2} &:=  1, &\text{for }& \ i=N,2N,\dots ,N^2  \\
b_{i3} &:=  1, &\text{for }& \ i=N(N-1)+1,N(N-1)+2,\dots ,N^2  \\
b_{i4} &:= 1,  &\text{for }& \ i=1,N+1,\dots,N(N-1)+1.                              
\end{align*}
In our example, we only use the first input, i.e., $u_2 \equiv u_3 \equiv u_4 \equiv 0$, $N=3$, define the output to be $y(t) = \sum_{i=1}^{N^2}T_i(t)$ and identify this model through the MATLAB System Identification Toolbox. In particular, we use the \emph{N4SID} algorithm to identify a continuous-time state-space system without noise model, where the input and noisy output is shown in \cref{fig:io_id}. Since the identified model is not externally positive, we then use \cref{alg:alt_cor} to find a near externally positive approximation of the identified model. The impulse response of the difference to the true model is shown in \cref{fig:imp_ex_pos} for both identified models. %

\subsection{State-feedback controller design}
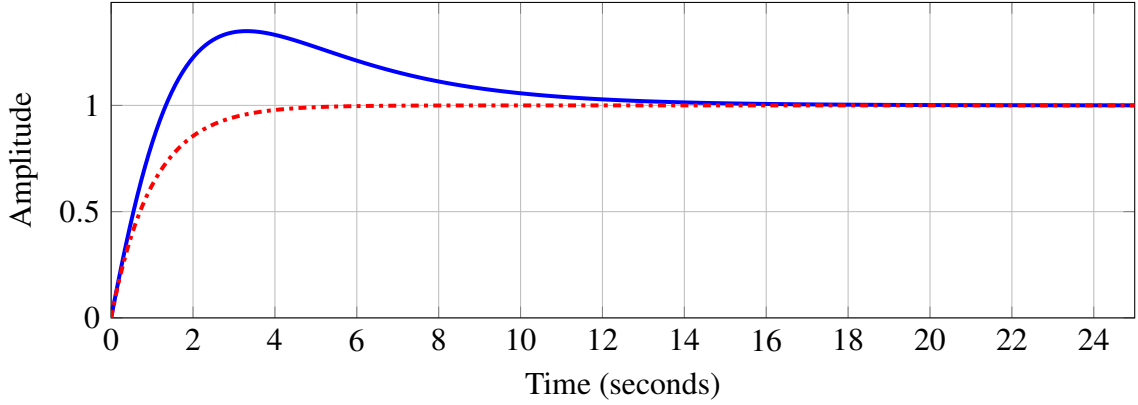
\begin{figure}
	\centering
	\begin{tikzpicture}
	\begin{groupplot}[group style={group name=my plots, group size=1 by 1}]
	
	\nextgroupplot[xlabel=Time (seconds),
	ylabel= Amplitude,
	xmin = 0,
	xmax = 25,
	grid = both,
	height = 0.35 \linewidth,
	width = .92 \linewidth,
	ymin = 0,
	scaled ticks=false, 
	]
	
	\addplot[color = FigColor1] file{step_cont_no_pos.txt};
	\label{line:cont_no_pos};
	
	\addplot[color = FigColor2,mark options = {solid}, dashdotted] file{step_cont_pos.txt};
	\label{line:cont_pos};
	\end{groupplot}
	\end{tikzpicture}
	\caption{The step responses of state-feedback controlled closed loop system \cref{eq:sys_closed}: \ref{line:cont_pos} controller design with \cref{alg:alt_cor,eq:stab_cont_mod} such that \cref{eq:sys_closed} is externally positive for $\gamma = 1$, \ref{line:cont_no_pos} controller design with \cref{eq:stab_cont_mod}, only.  \label{fig:state_feedback}}
\end{figure}

Let us next use \cref{alg:alt_cor}, to design a state-feedback control law with integral action 
\begin{align*}
u(t) &= L x(t) + l_i x_i(t) = \Delta A x_e(t) \\
\dot{x}_i(t) &= r(t)-y(t)
\end{align*}
that provides closed-loop external positivity for a two-compartment systems
\begin{align*}
\dot{x}_1(t) &= -x_1(t) + u(t)\\
\dot{x}_2(t) &= x_1(t) - x_2(t) + u(t)\\
y(t) &= x_2(t).
\end{align*}
The closed-loop dynamics are then given by
\begin{align}
\dot{x}_e (t) &= (A+F \Delta A) x_e(t) + B r(t)\\
y(t)  & = Cx_e(t) \label{eq:sys_closed}
\end{align}
where 
\begin{align*}
A = \begin{pmatrix}
-1 & 0 & 0 \\
1 & -1 & 0 \\
0 & -1 & 0 \\
\end{pmatrix}, \quad F = \begin{pmatrix}
1\\1\\0
\end{pmatrix}, \quad B = \begin{pmatrix}
1\\1\\1\\
\end{pmatrix}, \quad C = \begin{pmatrix}
0\\1\\0
\end{pmatrix}^\transp.
\end{align*}
Besides external positivity, we also require the real part of the closed loop poles to be smaller than 0.25, which we can achieve by \cref{rem:inertia} through modifying \cref{eq:stab_opt} to
\begin{equation}
\begin{aligned}
& \underset{Q, A_Q \in \Rnn}{\text{minimize}}  & & \|A_Q\| \\
& \text{subject to}  & & A Q+ Q A + A_Q^\transp F^\transp + F A_Q + 0.5 Q \preceq 0\\
& & & Q\succ 0.
\end{aligned} \label{eq:stab_cont_mod}
\end{equation}
The step responses of our controller versus a controller that is solely designed through \cref{eq:stab_cont_mod} is shown in \cref{fig:state_feedback}. 
Note that despite the simplicity of our setting, this is a fundamentally important problem. Control engineers often have to balance between the amount of overshoot and a fast response. However, if the references signal targets a limited capacity, any overshoot is inadmissible. 

\subsection{Model order reduction}
Finally, let us apply (positive) cone-balanced truncation to the above discretized heat equation model with $N=10$ and the use of the second and the fourth input, i.e., $u_1 \equiv u_3 \equiv 0$. The output $y$ is represented by the average temperature within five vertical stripes:
$$C = \blkdiag \left( \textbf{1}_{\frac{N^2}{5}}^T, \textbf{1}_{\frac{N^2}{5}}^T, \textbf{1}_{\frac{N^2}{5}}^T, \textbf{1}_{\frac{N^2}{5}}^T, \textbf{1}_{\frac{N^2}{5}}^T \right),$$
where $ \textbf{1}_{\frac{N^2}{5}} \in \mathbb{R}^{\frac{N^2}{5}}$ stands for the vector of all ones.

By \cite{grussler2012symmetry,grussler2012model,sootla2012scalable}, it is know that even a reduced model of order one often outperforms internally positivity preserving methods \cite{reis2009positivity,feng2010internal,li2011positivity,sootla2012scalable}. As the system has no symmetry as exploited in \cite{grussler2012model}, it suffices to compare our reduction approaches to balanced truncation (BT). 

\begin{figure}
	\centering
	\begin{tikzpicture}
	\begin{groupplot}[group style={group name=my plots, group size=1 by 1}]
	
	\nextgroupplot[xlabel=Reduced order $R$,
	ylabel=$\dfrac{\|G - G_R\|_\infty}{\|G\|_\infty}$,
	xmin = 1,
	xmax = 14,
	grid = both,
	ymode = log,
	height = 0.38 \linewidth,
	width = .92 \linewidth,
	xtick = {2,4,...,14},
	ytick = {1e-15,1e-10,1e-5,1e0},
	scaled ticks=false, 
	]
	
	\addplot[color = FigColor1,mark = *] file{heat_cone_error2_BT.txt};
	\label{line:heat_cone_error2_BT};
	
	\addplot[color = FigColor2,mark options = {solid}, mark = triangle*, dashdotted] file{heat_cone_error2_CBT.txt};
	\label{line:heat_cone_error2_CBT};
	
	\addplot[color = FigColor3,dashed, mark=diamond*,mark options = {solid}] file{heat_cone_error2_PCBT.txt};
	\label{line:heat_cone_error2_PCBT};
	\end{groupplot}
	\end{tikzpicture}
	\caption{Normalized $H_\infty$-error for discretized heat equation: \hfill \break 
		\ref{line:heat_cone_error2_BT} BT: standard balanced truncation \hfill \break \ref{line:heat_cone_error2_CBT} CBT: cone preserving balanced truncation, $\gamma= 0.2961$ \hfill \break \ref{line:heat_cone_error2_PCBT} PCBT: positivity preserving CBT, $\gamma = 0.2961$ \label{fig:error_heat2}}
\end{figure}
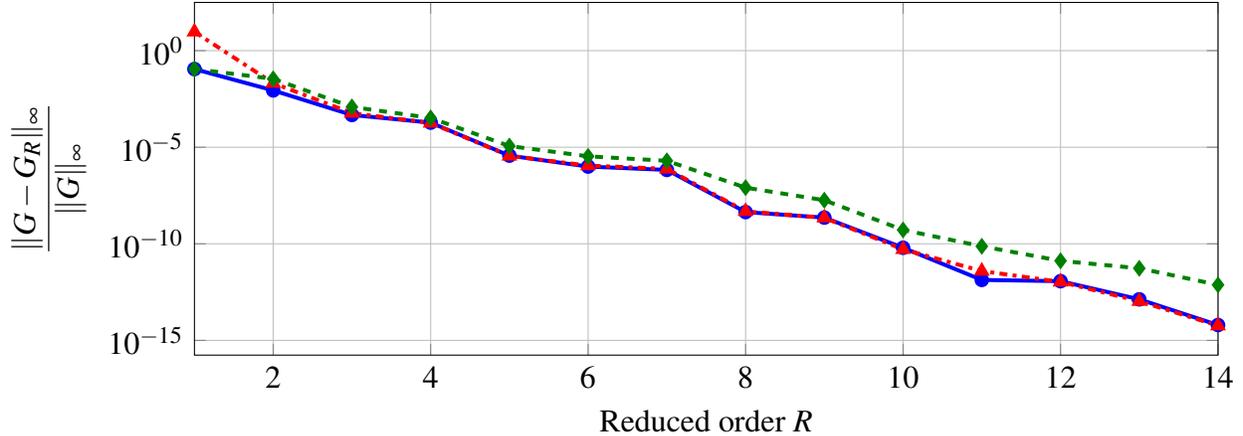

Our comparison starts from a minimal realization, which can be considered a \linebreak pre-reduction. For comparability of our results, we use the minimizing objective $\trace(Q+\tau C^TC)$ in case of positive cone balanced truncation (PCBT), which interestingly adds to an improved quality. For cone balanced truncation (CBT) we use the same $\gamma$ as determined by PCBT and $K$ is given by the equation$$A^TK + KA + 2\gamma K = -C^TC.$$ 
The normalized errors shown in~\cref{fig:error_heat2}, indicating that PCBT and CBT perform fairly close to BT.

%% file: sec_concl.txt
\section{Conclusion}
In this work, we have derived an external positivity certificate that is based on seeking invariant second-order cones. Our certificate has, in contrast to seeking an internally positive realization through invariant polyhedral cones, the advantage that it is tractable through semi-definite programming, has fixed computational cost and allows us to certify external positivity, where no other invariant cone would work. In particular, since our certificate is compatible with the SDP literature in control, we were able to exemplify the potential of our certificate through: 
\begin{enumerate}
	\item Establishing a heuristic solution to the problems of avoiding over- and undershooting in state-feedback controller design as well as to include prior knowledge of external positivity in system identification. 
	\item Modification of balanced truncation to preserve a dominant real pole.
	\item Modification of balanced truncation to preserve external positivity when our certificate applies.
\end{enumerate}
As indicated by our numerical examples, the established heuristic provides an interesting way of dealing with overshooting in the design of state-feedback controllers with integral action and helps to improve the quality of identified models. Further, our modified balanced truncation methods yield approximations that are qualitatively close to traditional balanced truncation. Thus, suggesting that these approaches only impose a mild conservatism in contrast to methods that preserve internal positivity \cite{reis2009positivity,sootla2012scalable}. However, it remains to understand how one can systematically choose a second-order cone that gives small truncation errors/error bounds.

{We have further provided justification why our certificate is likely to be less restrictive than existing ones. Nonetheless, our certificate is still only a sufficient test as we were able to construct systems that do not fulfil its requirements. In particular, we have constructed systems whose external positivity cannot be certified via second-order cones, polyhedral cones or any of the other mentioned certificates.} Thus, showing the need for more general cones, whose exponentially invariant certificate remains tractable. %

%% file: sec_appendix.txt
\section*{Appendix}
\appendix

\section{Proof to \cref{lem:dual}}
\label{proof:lem:dual}
	We start with the equivalence of \Cref{item:hyper_p,item:dual_cone}. For $\mathcal{K}_{K,c} = \mathcal{K}_{K,p}$ to hold true, $c$ must fulfill \cref{eq:strict_sep_equiv}, which by the definition of the dual cone \cref{eq:dual_cone} is the case if and only if $c \in\inter{\mathcal{K}_{K,p}^\ast}$. 	In order to see the set equality in \Cref{item:dual_cone}, note that there exists a $T \in \mathbb{R}^{n \times n}$ such that $T^\transp K T = K_n := \diag(1,\dots,1,-1)$ and $T^{-1}$ maps $\mathcal{K}_{K,c}$ onto the self-dual cone $\mathcal{K}_{K_n,e_n}$ \cite[Example~2.25]{boyd2004convex}, where $e_n$ is the n-th canonical unit vector. Hence, 
	\begin{equation*}
	(T^{-1} \mathcal{K}_{K,p}) = \mathcal{K}_{K_n,e_n} = (T^{-1} \mathcal{K}_{K,p})^\ast  = T^{\transp} \mathcal{K}_{K,p}^\ast,
	\end{equation*}
	and thus $\mathcal{K}_{K,p}^\ast = \mathcal{K}_{TK_nT^{\transp},T e_n} =   \mathcal{K}_{K^{-1},T e_n}$.  Then, as before, all normals to strictly separating hyperplanes of $\mathcal{K}_{K^{-1}}$ are given by $\inter{\mathcal{K}_{K,p}} = \inter{(\mathcal{K}_{K,p}^\ast)^\ast}$. 
	
	Then \Cref{item:dual_normal} just makes explicit \Cref{item:dual_cone}. Further, $c \in \inter{\mathcal{K}_{K,p}^\ast}$ if and only if there exists $\tau > 0$ such that $$\forall x \in \mathcal{K}_{K,p} \setminus \{0\} : x^\transp Kx + \tau x^\transp cc^\transp x > 0,$$ which is equivalent to
	$K+ \tau cc^\transp  \succ 0.$ Finally, by the inertia of $K$, the only admissible $\tau \in \mathds{R}$ to fulfil this inequality is a positive one.

\section{Proof to \cref{lem:inv_ellip}}
\label{proof:lem:inv_ellip}
The first part in \cref{eq:inv_ellip} is equivalent to $\mathcal{K}_K$ being $e^{At}$-invariant \cite{stern1991exponential}. To see the second part, note that by assumption $\mathcal{K}_{K,c}$ is proper, which by \cref{eq:strict_sep_equiv} means that there exists an $x \in \mathcal{K}_{K,c}: \; c^\transp x > 0$. Thus, the last item in \cref{lem:dual} (with $c = p$) applies.

\section{Proof to \cref{prop:necc_second}}
\label{proof:prop:necc_second}
	Without loss of generality, let $$A = \begin{pmatrix}
	0 & 0 & 0\\
	0 & 0 & \beta \\
	0 & -\beta & 0
	\end{pmatrix}$$
	Then for all $b \in \mathbb{R}^3$ with $b_1 \neq 0$, the set $\{e^{At}b: t \geq 0 \}$ is an ellipse, which implies that $\mathcal{K}_r(A,b) = \mathcal{K}_{K_b, c_b}$ with $K_b = \diag(-\frac{b_2^2 + b_3^2}{b_1^2},1,1)$ and $c_b = \sign(b_1) (1,0,0)^\transp$. Hence, any $e^{At}$-invariant proper convex cone $\mathcal{K}$ can be written as
	\begin{equation}
	\mathcal{K} = \cone(\bigcup_{b \in \mathcal{K}}  \mathcal{K}_{K_b, c_b}) =: \mathcal{K}_{K_{b_{\max}}, c_{b_{\max}}},
	\end{equation}
where $b_{\max} = \argmax_{b \in \mathcal{K}} \frac{b_2^2 + b_3^2}{b_1^2}$.

\section{Proof to \cref{thm:ex_pos_test}}
\label{proof:thm:ex_pos_test}
We begin by noticing that \cref{eq:thmtest:inertia_con}, this is, $\iota(K) = (n-1,0,1)$, implies the existence of a $p \in \mathds{R}^n$ such that $\mathcal{K}_{K,p}$ is a proper cone. Since \cref{eq:thmtest:b_in_K} is equivalent to $B_{(:,1)} \in \mathcal{K}_K$, we can assume that $p$ is chosen such that $B_{(:,1)} \in \mathcal{K}_{K,p}$. Therefore, if we can show that $C_{(i,:)} B_{(:,1)} > 0$ for all $i$, then \cref{eq:thmtest:c_dual} and the last item in \cref{lem:dual} allow us to conclude that $\mathcal{K}_{K,p} = \mathcal{K}_{C_{(i,:)}}$ for all $i$. To see this, note that \cref{eq:thmtest:c_dual} can only be fulfilled then if $\tau_i > 0$, which by $B_{(:,j)}^\transp (K+\tau_i C_{(i,:)}^\transp C_{(i,:)} ) B_{(:,j)} > 0$ and \cref{eq:thmtest:pos_init} yields that $C_{(i,:)}B_{(:,j)} > 0$ for all $j$ and $i$. In particular, this implies that $B_{(:,j)} \in \mathcal{K}_{C_{(1,:)}}$ for all $j$. Finally, by the second item in \cref{lem:dual} we have that $C_{(i,:)} \in \inter{\mathcal{K}_{C_{(1,:)}}^\ast}$ for all $i$ and $\mathcal{K}_{C_{(1,:)}}$ is $e^{At}$-invariant by \cref{lem:inv_ellip}. Thus proving that $(A,B,C)$ is $\mathcal{K}_{C_{(1,:)}}$-positive with strictly positive impulse response, which by \cref{prop:ext_pos_equiv} is equivalent to external positivity. 

\section{Proof to \cref{prop:sufficientonly}}
\label{proof:prop:sufficientonly}
	Without loss of generality, let $$A = \begin{pmatrix}
	\alpha & \beta & 0\\
	-\beta & \alpha & 0\\
	0 & 0 & 0
	\end{pmatrix} =: \blkdiag(A_1,0)$$
	where $\alpha < 0$ and $\beta \neq 0$. %
	Further, let $B = (1,0,1)^\transp$ and $x(t) := e^{At}B$. %
	Since \linebreak $\mathcal{S}:=\{(x_1(t),x_2(t)): t\geq 0\} \subset \{x \in \mathbb{R}^2: x_1^2+x_2^2 = 1 \}$ is not a closed contour, there exists a tangent hyperplane $\mathcal{T}^2 := \{x \in \mathbb{R}^2: c^T x \geq c_1\}$ to $\mathcal{C}^2 := \cl{\conv(\{(x_1(t),x_2(t)): t\geq 0\})}$ such that 
	\begin{enumerate}
		\item $\mathcal{C}^2 \subset \mathcal{T}^2$
		\item $\exists t^\opts > 0: c_1x_1(t^\opts)+c_2x_2(t^\opts) = c_1$
	\end{enumerate}
	and therefore $\mathcal{T} := \{x: Cx \geq 0 \}$ with $C := (c_1,c_2,-c_1)$ is a tangent hyperplane to $\mathcal{K}_r(A,B) = \cl{\cone( \{1\} \times \mathcal{C}^2)}$. Thus, $Ce^{At}B \geq 0$ for all $t\geq0$ and $CB = Ce^{At^\opts}B = 0$. In particular, for all $\tilde{B}_2 \notin \cl{\conv(\mathcal{S})}$ there exists $\tilde{t}\geq 0$ such that $C e^{A\tilde{t}}(\tilde{B}_2,1)^\transp < 0$ and thus $\mathcal{K}_r(A,B)$ is the only $(A,B)$-invariant cone that is contained in $\mathcal{T}$. Moreover, since $\conv(\mathcal{C}^2)$ is neither a polygon nor an ellipse, $\mathcal{K}_r(A,B)$ can neither be polyhedral nor second-order.
	
	Finally, note that for arbitrary $\varepsilon > 0$ and $\Delta C := (0,0,\varepsilon c_1)$, it holds that $(C+\Delta C)e^{At}B > 0$ for all $\forall t \geq 0$. Assume that for all $\varepsilon > 0$, there exists an $(A,B)$-invariant proper second order cone $\mathcal{K}_{K,p} \subset  \{x: (C+\Delta C) x \geq 0\}$. Then $\mathcal{K}_{K,p} \cap \{x: x_3 = 1\} =  \mathcal{E} := \{x: (x-k)^\transp P(x-k) \leq 1 \}$ for some $P \succ 0$ and $k \in \mathbb{R}^2$ with $\conv(\{x(t^\opts),(1,0)^\transp \}) \subset \mathcal{E} \subset \{x: c^\transp x \geq (1-\varepsilon) c_1\}$. However, as $\varepsilon \to 0$, this requires that either $\lambda_1(P) \to \infty$ or $\lambda_2(P) \to 0$. Thus the area of $\mathcal{E}$ can be made arbitrarily small or large, which either contradicts that $\mathcal{K}_{K,p}$ is $(A,B)$-invariant or $\mathcal{K}_{K,p} \subset  \{x: (C+\Delta C) x \geq 0\}$.

\section{Proof to \cref{thm:coneBT}}
\label{proof:thm:coneBT}

	The first part and the error bound follows as for generalized balanced truncation \cite{beck1996model,sandberg2004balanced}. \Cref{item:thmCBT:ABinv} follows by \begin{align*}
     \tilde{B}_{(:,j)}^\transp K \tilde{B}_{(:,j)} \geq  \tilde{B}_{(1:R,j)}^\transp K_{(1:R,1:R)} \tilde{B}_{(1:R,j)} \text{ for all } j \\
	\end{align*}
	which implies that if \cref{eq:thmtest:inv_K,eq:thmtest:b_in_K,eq:thmtest:inertia_con} are fulfilled for $(\tilde{A},B,C,D)$, $K$ and some $\gamma$, then the same applies to $(\tilde{A}_{(1:R:1:R)},\tilde{B}_{(1:R,:)}$. If additionally \cref{eq:thmtest:c_dual,eq:thmtest:pos_init} hold, then \cref{lem:dual} yields that
	\begin{align*} &(\tilde{B}_{(1:R,j)}^\transp \ 0)^\transp \in \mathcal{K}_{K,e_1} = \mathcal{K}_{K,C(i,:)} \text{ for all }i,j\\
	 &0 > \tilde{C}_{(i,:)} K^{-1}\tilde{C}_{(i,:)}^\transp \geq \tilde{C}_{(i,1:R)} K_{(1:R,1:R)} ^{-1}\tilde{C}_{(i,1:R)}^\transp \text{ for all } i,
	\end{align*}
	which shows that \cref{eq:thmtest:c_dual,eq:thmtest:pos_init} also hold for $(\tilde{A}_{(1:R:1:R)},\tilde{B}_{(1:R,:)},\tilde{C}_{(:,1:R)},D)$ and therefore external positivity is certified by \cref{thm:ex_pos_test}. Finally, \Cref{item:thmCBT:gamma} is obvious.

\section{Proof to \cref{thm:balancing}}
\label{proof:thm:balancing}

	Let $P$ and $K$ be as assumed. Then, we define $L := U \Sigma_P^{\frac{1}{2}}$ and $T := LV\Sigma^{-\frac{1}{2}}$ through the singular value decompositions of $P = U \Sigma_P U^{T}$ and {{eigenvalue decomposition of $L^{T}KL = V \bar{\Sigma} V^{T}$}} such that $\tilde{P} := T^{-1}PT^{-T}$ and $\tilde{K} := T^{T}QT$ fulfil
	\begin{align*}
	\tilde{P} & = \Sigma^{\frac{1}{2}}V^{T}L^{-1}LL^{T}L^{-T}V\Sigma^{\frac{1}{2}} = \Sigma, \\
	|\tilde{K}| &= |\Sigma^{-\frac{1}{2}}V^{T}L^{T}QLV\Sigma^{-\frac{1}{2}}| = \Sigma,
	\end{align*}
	with $\Sigma = \blkdiag(\sigma_1 I_{l_1},\dots,\sigma_s I_{l_s})$,  $\sigma_1 > \dots > \sigma_s > 0$, $l_1 + \dots + l_s = n$ and {{$\Sigma^2 = |\bar{\Sigma}|.$}} In particular, $\tilde{P}$ and $\tilde{K}$ are equal up to a sign-change on one of the diagonal entries by \cite[Theorem~4.5.8]{horn2012matrix}.
	
	Let us now verify that $\trace(\tilde{K}) < 0$, implying that the sign-change occurs at $\sigma_1$ with $l_1 = 1$. W.l.o.g, we assume that $P = I_n$ and $|K| = \Sigma^2$, i.e., 
	\begin{align}
	& \tilde{A}^T K + K \tilde{A} + 2\gamma K \preceq 0, \label{eq:I1}\\
	& \trace(NK) \leq 0\\
	& \tilde{A} + \tilde{A}^T =  -N.
	\end{align}
	Since, substituting $\tilde{A} = -N-\tilde{A}^T$ in~\cref{eq:I1} gives
	{\begin{align}
		-(N + \tilde{A}) K - K(N + \tilde{A}^T) -2 \gamma K &\preceq -4 \gamma K. \label{eq:tr}
		\end{align}}
	it follows by taking $\trace$ over $\cref{eq:tr}$ and the properties
	\begin{itemize}
		\item $\trace(NK) = \trace(KN)$
		\item $\trace(\tilde{A}K + K\tilde{A}^T + 2\gamma K) = \trace(\tilde{A}^TK +K\tilde{A} + 2\gamma K) \leq 0$
	\end{itemize}
	that
	\begin{align}
	2 \gamma \trace(K) \leq \trace(NK) \leq 0
	\end{align}
 The inertia of $K$ and the assumption that $\sigma_1 > \cdots > \sigma_s > 0$ imply then that the largest magnitude in $K$ is negative. The remaining part follows through state-space transformation.

\section{Proof to \cref{prop:delta}}
\label{proof:prop:delta}

	Let $(\tilde{A},\tilde{C})$ and $\tilde{N} = LL^\transp$ be as in the assumptions. Since \linebreak $\trace(N \tilde{K}) = \sum_{j} L_{(:,j)}^\transp \tilde{K} L_{(:,j)}$, we assume w.l.o.g. that $L_{(:,1)}^\transp \tilde{K} L_{(:,1)} < 0$. Thus, by \cref{lem:dual} and \cref{lem:inv_ellip} there exists a sufficiently large ${\varepsilon} > 0$ such that
	\begin{subequations}
		\begin{align}
		& \tilde{A}^\transp\tilde{K} + \tilde{K}\tilde{A} + 2\gamma \tilde{K} \preceq 0, \label{eq:1}\\
		& \tilde{A}\tilde{P} + \tilde{P}\tilde{A}^\transp \preceq -L_{(:,1)} L_{(:,1)}^\transp , \label{eq:2}\\
		& \tilde{K}^{-1} + {\varepsilon} L_{(:,1)} L_{(:,1)}^\transp  \succ 0,\label{eq:3}\\
		& 2\gamma{\varepsilon} p_{11} - p_{11}^{-1} > 0, \label{eq:4}
		\end{align}
	\end{subequations}
	Multiplying \cref{eq:1} with $\tilde{K}^{-1}$ from the right and the left yields 
	\begin{equation}
	\tilde{A}\tilde{K}^{-1} + \tilde{K}^{-1}\tilde{A}^\transp +2\gamma \tilde{K}^{-1} \preceq 0 \label{Qlmi}
	\end{equation}
	and multiplying \cref{eq:2} by $2\gamma {\varepsilon}$ gives 
	\begin{equation}
	2\gamma {\varepsilon} \tilde{A} \tilde{P} + 2\gamma {\varepsilon} \tilde{P} \tilde{A}^\transp  + 2 \gamma {\varepsilon} L_{(:,1)} L_{(:,1)}^\transp \preceq 0, \label{Pgram}
	\end{equation}
	Adding up \cref{Qlmi} and \cref{Pgram} results in
	\begin{equation*}
	\tilde{A}\Delta^{-1} + \Delta^{-1}\tilde{A}^\transp + 2\gamma \left(\tilde{K}^{-1} + {\varepsilon} L_{(:,1)} L_{(:,1)}^\transp  \right) \preceq 0
	\end{equation*}
	with $\Delta := (2\gamma {\varepsilon} \tilde{P} + \tilde{K}^{-1})^{-1} \succ 0$. Finally, a proper scaling of $\Delta$ gives a diagonal solution to
	\begin{equation}
	\tilde{A}^\transp\Delta + \Delta\tilde{A} \preceq -\tilde{C}^\transp\tilde{C} \label{eq:D}.
	\end{equation}
	The last implication follows by \cref{thm:balancing}.